\documentclass{amsart}

\usepackage{amssymb}
\usepackage{amsmath}
\usepackage[all]{xy}

\newtheorem{theorem}[subsection]{Theorem}
\newtheorem*{ntheorem}{Theorem}
\newtheorem{cor}[subsection]{Corollary}

\newtheorem{lemma}[subsection]{Lemma}

\newtheorem{proposition}[subsection]{Proposition}
\newtheorem*{nproposition}{Proposition}

\theoremstyle{remark}
\newtheorem{exa}[subsection]{Example}
\newtheorem{rem}[subsection]{Remark}

\theoremstyle{definition}

\DeclareMathOperator{\Int}{Int}

\DeclareMathOperator{\codim}{codim}
\DeclareMathOperator{\BOX}{Box}
\DeclareMathOperator{\Hom}{Hom}

\DeclareMathOperator{\vol}{vol}

\DeclareMathOperator{\orb}{orb}

\DeclareMathOperator{\com}{c}
\DeclareMathOperator{\rk}{rk}

\newcommand{\bSigma}{\mbox{\boldmath$\Sigma$}}
\newcommand{\btau}{\mbox{\boldmath$\tau$}}

\newcommand{\btriangle}{\mbox{\boldmath$\triangle$}}

\begin{document}

\title{Weighted Ehrhart Theory and Orbifold Cohomology}

\author{A. Stapledon}

\address{Department of Mathematics, University of Michigan, Ann Arbor, MI 48109, USA}
\email{astapldn@umich.edu}

\begin{abstract}
We introduce the notion of a weighted $\delta$-vector of a lattice polytope.
Although the definition is motivated by motivic integration, we study weighted $\delta$-vectors from a combinatorial perspective. We present a version of Ehrhart Reciprocity and
prove a change of variables formula. We deduce a new geometric  interpretation of the coefficients of the Ehrhart $\delta$-vector.
More specifically, they are sums of dimensions of orbifold cohomology groups of a toric stack.
\end{abstract}

\maketitle

\section{Introduction}

Let  $P$ be a $d$-dimensional lattice polytope in a lattice $N$ of rank $d$.
For each positive integer $m$, let $f_{P}(m)$ be the number
of lattice points in $mP$. Then $f_{P}(m)$ is a polynomial in $m$ of degree $d$, called the
\emph{Ehrhart polynomial} of $P$ \cite{EhrLinearI,EhrLinearII}.
The generating series of the Ehrhart polynomial can be written in the form
\begin{equation*}
\sum_{m \geq 0} f_{P}(m)t^{m} = \delta_{P}(t)/ (1 - t)^{d + 1},
\end{equation*}
where $\delta_{P}(t)$ is a polynomial of degree less than or equal to $d$ with non-negative integer coefficients
\cite{StaDecompositions}. With a slight abuse of terminology,
$\delta_{P}(t)$ is called the \emph{(Ehrhart) $\delta$-vector} of $P$.
We will write
$
\delta_{P}(t) = \delta_{d}t^{d} + \delta_{d - 1}t^{d-1} + \cdots + \delta_{0}$.

We present a new geometric interpretation of the coefficients $\delta_{i}$ of the Ehrhart $\delta$-vector.
When $P$ is reflexive, Batyrev and Dais \cite{BDStrong} showed that $\delta_{i}$ is the $2i^{\textrm{th}}$ stringy Betti number of a toric variety $X$ associated to $P$.
Furthermore, in this case, results of Yasuda \cite{YasMotivic} imply that $\delta_{i}$ is equal to the dimension of the $2i^{\textrm{th}}$ orbifold cohomology group of the canonical orbifold associated to $X$. This interpretation of the Ehrhart $\delta$-polynomial of $P$ was used by Musta\c t\v a and Payne in \cite{MPEhrhart} and Karu in \cite{KarEhrhart}.
We generalise this result to any lattice polytope $P$ by showing that $\delta_{i}$ is a sum of dimensions of orbifold cohomology groups of a toric stack (Theorem \ref{boo}).

In order to establish this result, we introduce certain refinements of the Ehrhart polynomial. Fix a lattice point $\alpha$ in $P$.
After translating, we may assume that $\alpha$ is the origin. Denote the union of the facets of $P$ not containing the origin by $\partial P_{0}$.
Consider the fan $\triangle$ over the faces of $\partial P_{0}$, with support $|\triangle|$. 
Let $\psi: |\triangle| \rightarrow \mathbb{R}$ be the piecewise $\mathbb{Q}$-linear function with respect to $\triangle$
satisfying $\psi(v) = 1$ for all $v$ in $\partial P_{0}$. We fix a simplicial fan $\Sigma$ refining $\triangle$ and with the same rays as $\triangle$ (see, for example, \cite{OPLinear}).
Let $\lambda: |\Sigma| \rightarrow \mathbb{R}$ be an arbitrary piecewise $\mathbb{Q}$-linear function with respect to $\Sigma$ satisfying 
$\lambda(v) > -1$ for all $v$ in $\partial P_{0}$.
Consider the function
\begin{equation*}
w_{\lambda}: |\Sigma| \cap N \rightarrow \mathbb{Q}
\end{equation*}
\[
w_{\lambda}(v) = \psi(v) - \lceil \psi(v) \rceil + \lambda(v).
\]
Note that when $\lambda \equiv 0$, the function $w_{0}$ is determined by the pair $(P, \alpha)$.
Think of $w_{\lambda}$ as assigning a weight to every lattice point in $|\Sigma|$.
For every rational number $k$ and for every non-negative integer $m$,
denote by $f_{k}^{\lambda}(m)$ the number of lattice points of weight $k$ in $mP$.
Note that the Ehrhart polynomial of $P$ can be recovered as $f_{P}(m) = \sum_{k \in \mathbb{Q}} f_{k}^{\lambda}(m)$.
For every rational number $k$, consider the power series 
\begin{equation*}
\delta^{\lambda}_{k}(t) := (1 - t)^{d + 1} \sum_{m \geq 0} f_{k}^{\lambda}(m)t^{m}.
\end{equation*}
If $\lambda(v) \geq 0$  for all $v$ in $|\Sigma| \cap N$, then $\delta^{\lambda}_{k}(t)$ is a polynomial in $t$ with integer coefficients 
(Corollary \ref{Coral}).
The $\delta$-vector of $P$ decomposes as $\delta_{P}(t) = \sum_{k \in \mathbb{Q}} \delta^{\lambda}_{k}(t)$.
We define the \emph{weighted $\delta$-vector} of $P$ to be
\begin{equation*}
\delta^{\lambda}(t) := \sum_{k \in \mathbb{Q}} \delta^{\lambda}_{k}(t) t^{k}.
\end{equation*}
In Section \ref{weight}, we 
verify that $\delta^{\lambda}(t)$ lies in  $\mathbb{Z}[[t^{1/N}]]$, for some positive integer $N$.
We discuss the geometric interpretation of the weighted $\delta$-vector in terms of motivic integration
in \cite{YoWeightII}.

We emphasise that the weighted $\delta$-vector of $P$ is interesting from a purely combinatorial perspective.
We present a change of variables formula in Section \ref{weight}, comparing weighted $\delta$-vectors on different polytopes (Proposition \ref{change}).
The weighted $\delta$-vector has the following symmetry property (Corollary \ref{beautiful}).
\begin{nproposition}
If the origin lies in the interior of $P$,
then the weighted $\delta$-vector $\delta^{\lambda}(t)$ 
satisfies
\begin{equation*}
\delta^{\lambda}(t) = t^{d}\delta^{\lambda}(t^{-1}).
\end{equation*}
\end{nproposition}

The most important example is when
$\lambda \equiv 0$, since the Ehrhart $\delta$-vector $\delta_{P}(t)$ can be easily recovered from
the weighted $\delta$-vector $\delta^{0}(t)$ (see (\ref{bling})).
More specifically, $\delta^{0}(t)$ is a polynomial of degree $d$ with rational powers and non-negative integer coefficients (see (\ref{equat}))
and the coefficient of $t^{i}$ in  $\delta_{P}(t)$ is the sum of the coefficients of
$t^{j}$ in $\delta^{0}(t)$ for $i - 1 < j \leq i$.
In Section \ref{Ehr}, we consider this case and deduce the following weighted version of Ehrhart Reciprocity (Theorem \ref{propB}).



\begin{ntheorem}[Weighted Ehrhart Reciprocity]
Suppose the origin lies in the interior of $P$.
For every rational number $- 1 < k \leq 0$, $f^{0}_{k}(m)$ is either identically zero or
a polynomial of degree $d$ in
$\mathbb{Q}[t]$ with positive leading coefficient.
For any positive integer $m$,
\begin{displaymath}
f^{0}_{k}(-m) = \left\{ \begin{array}{ll}
(-1)^{d}f^{0}_{-1-k}(m) & \textrm{ if } -1 < k < 0 \\
(-1)^{d}f^{0}_{k}(m - 1) & \textrm{ if } k = 0.
\end{array} \right.
\end{displaymath}
\end{ntheorem}

We show that Ehrhart Reciprocity for lattice polytopes is an immediate consequence (Corollary \ref{Ehrt}) as well as a result of Hibi (Corollary \ref{gogo}).
In fact, Weighted Ehrhart Reciprocity implies Ehrhart Reciprocity for rational polytopes (Remark \ref{TMac}).




We now consider the geometric side of the story. Let $v_{1}, \ldots, v_{r}$ be the primitive integer vectors of $\Sigma$.
For $i = 1, \ldots, r$, there is a positive integer $a_{i}$ such that $a_{i}v_{i}$ lies in $\partial P_{0}$.
The data $\bSigma = (N, \Sigma, \{ a_{i}v_{i} \})$ is called a \emph{stacky fan}. We can associate to $\bSigma$ a Deligne-Mumford stack $\mathcal{X} = \mathcal{X}(\bSigma)$
with coarse moduli space
the toric variety $X = X(\Sigma)$ \cite{BCSOrbifold}.
The theory of orbifold cohomology, developed by Chen and Ruan \cite{CRNew, CROrbifold},
associates to $\mathcal{X}$ a finite-dimensional $\mathbb{Q}$-algebra
$H_{\orb}^{*}(\mathcal{X}, \mathbb{Q})$, graded by $\mathbb{Q}$.
In Section \ref{orbifold}, we use a result of Borisov, Chen and Smith (Proposition 4.7 \cite{BCSOrbifold})
to deduce our desired geometric interpretation of the coefficients of the $\delta$-vector of $P$ (Theorem \ref{boo}).

\begin{ntheorem}
The coefficient of $t^{j}$ in $\delta^{0}(t)$ is equal to $\dim_{\mathbb{Q}} H_{\orb}^{2j}(\mathcal{X}, \mathbb{Q})$. Moreover, the coefficient $\delta_{i}$ of $t^{i}$ in the $\delta$-vector
$\delta_{P}(t)$ is a sum of dimensions of orbifold cohomology groups,
\begin{equation*}
\delta_{i}  = \sum_{2i - 2< j \leq 2i} \dim_{\mathbb{Q}} H_{\orb}^{j}(\mathcal{X}, \mathbb{Q}).
\end{equation*}
\end{ntheorem}

Using the above theorem, 
we show that Weighted Ehrhart Reciprocity  is  a consequence of Poincar\'e duality for orbifold cohomology (Remark \ref{bias}).

Another corollary of the above theorem is that we may express the coefficients of the $\delta$-vector of $P$ as dimensions of orbifold cohomology groups of a $(d + 1)$-dimensional orbifold (Theorem \ref{chance}). More specifically, fix a lattice triangulation $\mathcal{T}$ of $P$ and let
$\sigma$ be the cone over $P \times \{ 1 \}$ in $(N \times \mathbb{Z})_{\mathbb{R}}$. If $\triangle$ denotes the simplicial fan refinement of $\sigma$ induced by $\mathcal{T}$, then we may consider the corresponding toric variety $Y = Y(\triangle)$, with its canonical orbifold structure.

\begin{ntheorem}
Let $P$ be a $d$-dimensional lattice polytope and let $\mathcal{T}$ be a lattice triangulation of $P$ corresponding to a $(d + 1)$-dimensional toric variety $Y$ as above. The Ehrhart $\delta$-vector of $P$ has the form
\[
\delta_{P}(t) = \sum_{i = 0}^{d} \dim_{\mathbb{Q}} H_{\orb}^{2i}(Y, \mathbb{Q}) t^{i}.
\]
\end{ntheorem}

In the final section, we give a third proof of Weighted Ehrhart Reciprocity, generalising the toric proof of Ehrhart Reciprocity in
\cite{FulIntroduction}. More specifically, we show that Weighted Ehrhart Reciprocity is a consequence of Serre Duality as well as some vanishing
theorems for ample divisors on toric varieties due to Musta\c t\v a \cite{MusVanishing}. We note that this proof applies under the assumptions of the introduction, but not under the more general assumptions in Theorem \ref{propB}.

We end the introduction with an example 
illustrating Weighted Ehrhart Reciprocity and showing how the weighted $\delta$-vector gives rise
to the Ehrhart $\delta$-vector.

\begin{exa}
Let $N = \mathbb{Z}^{2}$ and let $P$ be the lattice polytope with vertices $(1,0)$,$(0,2)$,
$(-1,2)$, $(-2,1)$, $(-2,0)$ and $(0,-1)$.
Since the origin lies in the interior of $P$, weighted Ehrhart Reciprocity holds and the weighted $\delta$-vector $\delta^{0}(t)$ is symmetric.
We can compute $\delta^{0}(t)$ using Lemma \ref{Kansas} and 
one can show that
\[ f^{0}_{0}(m) = |\partial P \cap N|m(m + 1)/2 + 1 \] \[ f^{0}_{k}(m) = (f^{0}_{k}(1) +  f^{0}_{-1 - k}(1)) m^{2}/2 + ( f^{0}_{k}(1) -  f^{0}_{-1 - k}(1)) m/2 \textrm{ for } k \neq 0. \]
Note that the lattice points of weight $0$ are precisely those lying on $\partial(mP)$ for some non-negative integer $m$. We have marked the
lattice points of non-zero weight in $2P$.
\[ \delta^{0}(t) = t^{2} + 2t^{3/2} + t^{4/3} + 4t + t^{2/3} + 2t^{1/2} + 1 \]
\[ \delta_{P}(t) = 4t^{2} + 7t + 1 \]

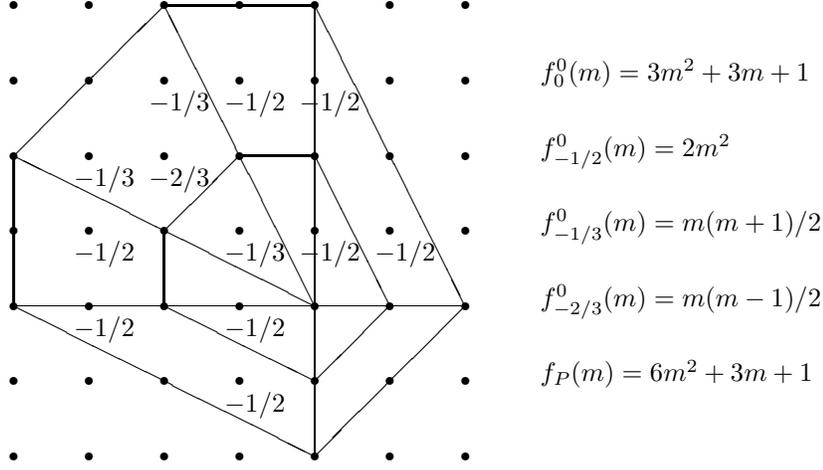
\setlength{\unitlength}{1cm}
\begin{picture}(12,6.5)

\put(7,5){$f^{0}_{0}(m) = 3m^{2} + 3m + 1$}
\put(7,4){$f^{0}_{-1/2}(m) = 2m^{2}$}
\put(7,3){$f^{0}_{-1/3}(m) = m(m + 1)/2$}
\put(7,2){$f^{0}_{-2/3}(m) = m(m - 1)/2$}
\put(7,1){$f_{P}(m) = 6m^{2} + 3m + 1$}

\put(0,0){\circle*{0.1}}
\put(0,1){\circle*{0.1}}
\put(0,2){\circle*{0.1}}
\put(0,3){\circle*{0.1}}
\put(0,4){\circle*{0.1}}
\put(0,5){\circle*{0.1}}
\put(0,6){\circle*{0.1}}
\put(1,0){\circle*{0.1}}
\put(1,1){\circle*{0.1}}
\put(1,2){\circle*{0.1}}
\put(1,3){\circle*{0.1}}
\put(1,4){\circle*{0.1}}
\put(1,5){\circle*{0.1}}
\put(1,6){\circle*{0.1}}
\put(2,0){\circle*{0.1}}
\put(2,1){\circle*{0.1}}
\put(2,2){\circle*{0.1}}
\put(2,3){\circle*{0.1}}
\put(2,4){\circle*{0.1}}
\put(2,5){\circle*{0.1}}
\put(2,6){\circle*{0.1}}
\put(3,0){\circle*{0.1}}
\put(3,1){\circle*{0.1}}
\put(3,2){\circle*{0.1}}
\put(3,3){\circle*{0.1}}
\put(3,4){\circle*{0.1}}
\put(3,5){\circle*{0.1}}
\put(3,6){\circle*{0.1}}
\put(4,0){\circle*{0.1}}
\put(4,1){\circle*{0.1}}
\put(4,2){\circle*{0.1}}
\put(4,3){\circle*{0.1}}
\put(4,4){\circle*{0.1}}
\put(4,5){\circle*{0.1}}
\put(4,6){\circle*{0.1}}
\put(5,0){\circle*{0.1}}
\put(5,1){\circle*{0.1}}
\put(5,2){\circle*{0.1}}
\put(5,3){\circle*{0.1}}
\put(5,4){\circle*{0.1}}
\put(5,5){\circle*{0.1}}
\put(5,6){\circle*{0.1}}
\put(6,0){\circle*{0.1}}
\put(6,1){\circle*{0.1}}
\put(6,2){\circle*{0.1}}
\put(6,3){\circle*{0.1}}
\put(6,4){\circle*{0.1}}
\put(6,5){\circle*{0.1}}
\put(6,6){\circle*{0.1}}

\put(3.8,2.6){$-1/2$}
\put(3.8,4.6){$-1/2$}
\put(4.8,2.6){$-1/2$}
\put(2.8,1.6){$-1/2$}
\put(0.8,1.6){$-1/2$}
\put(0.8,2.6){$-1/2$}
\put(0.8,3.6){$-1/3$}
\put(1.8,3.6){$-2/3$}
\put(1.8,4.6){$-1/3$}
\put(2.8,2.6){$-1/3$}
\put(2.8,4.6){$-1/2$}
\put(2.8,0.6){$-1/2$}


\linethickness{0.075mm}
\put(4,0){\line(0,1){6}}
\put(0,2){\line(1,0){6}}
\put(4,2){\line(-1,2){2}}
\put(4,2){\line(-2,1){4}}

\linethickness{0.2mm}
\put(3,4){\line(1,0){1}}
\put(2,2){\line(0,1){1}}
\put(2,6){\line(1,0){2}}
\put(0,2){\line(0,1){2}}
\put(4,4){\line(1,-2){1}}
\put(5,2){\line(-1,-1){1}}
\put(2,2){\line(2,-1){2}}
\put(2,3){\line(1,1){1}}
\put(4,6){\line(1,-2){2}}
\put(6,2){\line(-1,-1){2}}
\put(0,2){\line(2,-1){4}}
\put(0,4){\line(1,1){2}}

\end{picture}

\end{exa}

We give a brief outline of the paper and note that we will use a more general setup than that in
the introduction.
In Section \ref{weight}, we introduce the notion of a weighted $\delta$-vector and describe some of its properties. We specialise to the case when $\lambda \equiv 0$ and prove
weighted Ehrhart Reciprocity in Section \ref{Ehr}. In Section \ref{orbifold}, we consider orbifold cohomology and deduce our geometric interpretation of the coefficients of the $\delta$-vector of $P$. In Section \ref{toric}, we give a proof of Weighted Ehrhart Reciprocity using Serre Duality and vanishing theorems for ample divisors on toric varieties.

\medskip

The author would like to thank Mircea Musta\c t\v a for his constant help, encouragement and patience.
He is very grateful to Sam Payne for carefully reading several preliminary drafts and providing valuable feedback.
He would also like to thank Alexander Barvinok, Bill Fulton and Kevin Tucker for some useful discussions.
The author was supported by
Mircea Musta\c t\v a's Packard Fellowship and
by an Eleanor Sophia Wood
travelling scholarship from the University of Sydney.

\section{Weighted $\delta$-Vectors}\label{weight}

The goal of this section is to define weighted $\delta$-vectors and study their properties.
We will fix the following notation throughout the paper. Our setup will be slightly more general than that in the introduction (c.f. Remark \ref{impo}).
Let $N$ be a lattice of rank $d$ and set $N_{\mathbb{R}} = N \otimes_{\mathbb{Z}} \mathbb{R}$.
Let $\Sigma$ be a simplicial, rational, $d$-dimensional fan in $N_{\mathbb{R}}$.
We assume that the support $|\Sigma|$ of $\Sigma$ in $N_{\mathbb{R}}$ is convex. Recall that $\Sigma$ is \emph{complete} if $|\Sigma| = N_{\mathbb{R}}$.
Let $\rho_{1}, \ldots, \rho_{r}$ denote
the rays of $\Sigma$, with primitive integer generators $v_{1}, \ldots, v_{r}$ in $N$.
Fix elements $b_{1}, \ldots , b_{r}$ in $N$ such that $b_{i} = a_{i}v_{i}$ for some positive integer
$a_{i}$, for $i = 1, \ldots, r$.
The data $\bSigma = (N, \Sigma, \{ b_{i} \} )$ is called a \emph{stacky fan} \cite{BCSOrbifold}.
Let $\psi: |\Sigma| \rightarrow \mathbb{R}$ be the function that is $\mathbb{Q}$-linear on each cone of
$\Sigma$ and satisfies $\psi(b_{i}) = 1$ for $i = 1, \ldots ,r$. We define
\begin{equation*}
Q = Q_{\bSigma} = \{ v \in |\Sigma| \mid \psi(v) \leq 1 \}.
\end{equation*}
Observe that $Q$ need not be convex and
that the union of the facets in the boundary $\partial Q$ of $Q$ not containing $0$ is given by
\begin{equation*}
\{ v \in N_{\mathbb{R}} \mid \psi(v) = 1 \}.
\end{equation*}
Hence $\psi$ depends only on $Q$, not on the corresponding stacky fan.

\begin{rem}
A \emph{pure lattice complex} of dimension $d$ is a simplicial complex such that all maximal faces are lattice polytopes of dimension $d$ \cite{BMLattice}.
For any cone $\tau$ in $\Sigma$, $Q \cap \tau$ is a lattice simplex containing the origin as a vertex. It follows that $Q$ has the structure of a pure lattice complex.
Conversely, let $K$ be a pure lattice complex of dimension $d$ in $N_{\mathbb{R}}$ such that the smallest cone containing $K$ is convex. If every maximal simplex in $K$ contains the origin
as a vertex, then $K = Q_{\bSigma}$ for some stacky fan $\bSigma$ as above. More specifically,
let $\Sigma$ be the cone over the faces in $K$ not containing the origin and make an appropriate
choice of $\{b_{i}\}$.
\end{rem}

For each positive integer $m$, let $f_{Q}(m)$ be the number of lattice points in $mQ$. Then $f_{Q}(m)$ is a polynomial in $m$ of degree $d$, called the
\emph{Ehrhart polynomial} of $Q$ \cite{BMLattice}. We write
\begin{equation}\label{poly}
f_{Q}(m) = c_{d}m^{d} + c_{d - 1}m^{d - 1} + \cdots + c_{0}.
\end{equation}
The generating series of the Ehrhart polynomial can be written in the form
\begin{equation*}
\sum_{m \geq 0} f_{Q}(m)t^{m} = \delta_{Q}(t)/ (1 - t)^{d + 1},
\end{equation*}
where $\delta_{Q}(t)$ is a polynomial of degree less than or equal to $d$ with non-negative integer coefficients \cite{StaDecompositions}. With a slight abuse of terminology,
$\delta_{Q}(t)$ is called the \emph{(Ehrhart) $\delta$-vector} of $Q$.
We write
\begin{equation}\label{otation}
\delta_{Q}(t) = \delta_{d}t^{d} + \delta_{d - 1}t^{d-1} + \cdots + \delta_{0}.
\end{equation}

\begin{rem}\label{cough}
Consider a triple $(N, \triangle, \{b_{i}\})$ as above 
but without the assumption that $\triangle$ is simplicial.
Suppose there exists a piecewise $\mathbb{Q}$-linear function $\psi: |\triangle| \rightarrow \mathbb{R}$
satisfying $\psi(b_{i}) = 1$, and set $Q = \{ v \in |\triangle| \mid \psi(v) \leq 1 \}$.
There exists a simplicial fan $\Sigma$ refining $\triangle$ and with the same rays as $\triangle$ (see, for example, \cite{OPLinear}).
If $\bSigma$ denotes the stacky fan  $(N, \Sigma, \{ b_{i} \})$, then $Q = Q_{\bSigma}$.
\end{rem}

\begin{rem}\label{impo}
We consider the following important example. Let $P$ be a lattice polytope and let $\alpha$ be a lattice point in $P$. After translating, we may assume that $\alpha$ is the origin.
Let $\triangle$ be the fan over the faces of $P$ not containing the origin. As in Remark \ref{cough}, let $\Sigma$  be
a simplicial fan refining $\triangle$ and with the same rays as $\triangle$.
With the appropriate choice of  $\{b_{i}\}$, $P = Q_{\bSigma}$.  In the introduction, we state the main results of the paper in this context.
\end{rem}

Let $\lambda: |\Sigma| \rightarrow \mathbb{R}$ be an arbitrary piecewise $\mathbb{Q}$-linear function with respect to $\Sigma$.
We introduce the `weight function' $w_{\lambda}$ on $|\Sigma| \cap N$,
\begin{equation*}
w_{\lambda}: |\Sigma| \cap N \rightarrow \mathbb{Q}
\end{equation*}
\[
w_{\lambda}(v) = \psi(v) - \lceil \psi(v) \rceil + \lambda(v).
\]
Note that when $\lambda \equiv 0$,  the corresponding weight function $w_{0}$ is determined by $Q$.
Think of $w_{\lambda}$ as assigning a weight to every lattice point in $|\Sigma|$.
For every rational number $k$ and for every non-negative integer $m$,
denote by $f_{k}^{\lambda}(m)$ the number of lattice points of weight $k$ in $mQ$.
Note that the Ehrhart polynomial of $Q$ can be recovered as $f_{Q}(m) = \sum_{k \in \mathbb{Q}} f_{k}^{\lambda}(m)$.
For every rational number $k$, consider the power series 
\begin{equation*}
\delta^{\lambda}_{k}(t) := (1 - t)^{d + 1} \sum_{m \geq 0} f_{k}^{\lambda}(m)t^{m}.
\end{equation*}
If $\lambda(v) \geq 0$  for all $v$ in $|\Sigma|$, we will show that $\delta^{\lambda}_{k}(t)$ is a polynomial in $t$ with integer coefficients 
(Corollary \ref{Coral}).
The Ehrhart $\delta$-vector of $Q$ decomposes as $\delta_{Q}(t) = \sum_{k \in \mathbb{Q}} \delta^{\lambda}_{k}(t)$.
We define the \emph{weighted $\delta$-vector} of $Q$ by
\begin{equation*}
\delta^{\lambda}(s,t) := \sum_{k \in \mathbb{Q}} \delta^{\lambda}_{k}(t) s^{k}.
\end{equation*}
It follows from the definition that the weighted $\delta$-vector  can be written as
\[ \delta^{\lambda}(s,t) = (1 - t)^{d + 1} \sum_{m \geq 0} (\sum_{v \in mQ \cap N} s^{w_{\lambda}(v)}) t^{m},\]
and hence $\delta^{\lambda}(s,t) $ is a well-defined element of
$\mathbb{Z}[s^{q} \mid q \in \mathbb{Q}][[t]]$. Note that when $s = 1$, we recover the Ehrhart $\delta$-vector $\delta^{\lambda}(1,t) = \delta_{Q}(t)$.
Later, as in the introduction, we will consider the case when $\delta^{\lambda}(t,t)$ lies in $\mathbb{Z}[[t^{1/N}]]$ for some positive integer $N$.
Our first aim is to express $\delta^{\lambda}(s,t)$ as a rational function in $\mathbb{Q}(s^{1/N},t)$, for some positive integer $N$ (Proposition \ref{PropA}).



For each cone $\tau$ in $\Sigma$, let $\Sigma_{\tau}$ be the simplicial fan in $(N/N_{\tau})_{\mathbb{R}}$ with cones given by the projections
of the cones in $\Sigma$ containing $\tau$. If $\tau$ is not contained in the boundary of $|\Sigma|$, then $\Sigma_{\tau}$ is complete.
The $h$-vector of $\Sigma_{\tau}$ is given by
\[
h_{\tau}(t) := \sum_{\tau \subseteq \sigma} t^{\dim \sigma - \dim \tau}(1 - t)^{\codim \sigma}.
\]
We will sometimes write $h_{\Sigma}(t)$ for $h_{\{ 0 \}}(t)$.
We will use the following standard lemma. For a combinatorial proof, we refer the reader to  \cite{StaEnumerative} and Lemma 1.3 \cite{HibLower}. We provide a geometric proof, deducing the result as a corollary of Lemma \ref{AFL}, which is proved independently.

\begin{lemma}\label{Lefty}
For each cone $\tau$ in $\Sigma$, $h_{\tau}(t)$ is a polynomial of degree at most $\codim \tau$ with non-negative integer coefficients.
Suppose that $\tau$ is not contained in the boundary of $\Sigma$. Then 
$h_{\tau}(t) = t^{\codim \tau}h_{\tau}(t^{-1})$ and the coefficients of
$h_{\tau}(t)$ are positive integers.
\end{lemma}
\begin{proof}
It follows from the definition that $h_{\tau}(t)$ is a polynomial of degree at most $\codim \tau$.
Consider the simplicial toric variety $X = X(\Sigma_{\tau})$.
By Lemma \ref{AFL}, we can interpret the coefficient of $t^{i}$ in $h_{\tau}(t)$ as the dimension of the $2i^{\textrm{th}}$ cohomology group of $X$.
In particular, each coefficient is non-negative.
If $\tau$ is not contained in the boundary of $\Sigma$, then $X$ is complete. In this case,
$h_{\tau}(t) = t^{\codim \tau}h_{\tau}(t^{-1})$ follows from Poincar\'e duality on $X$ and
the coefficients of $h_{\tau}(t)$ are positive by Lemma \ref{AFL}.
\end{proof}

We define the \emph{weighted $h$-vector} $h_{\tau}^{\lambda}(s,t)$ of $\Sigma_{\tau}$,
\begin{equation}\label{pretty}
h_{\tau}^{\lambda}(s,t) :=  \sum_{\tau \subseteq \sigma}s^{\sum_{\rho_{i} \subseteq \sigma \backslash \tau} \lambda(b_{i})}  t^{\dim \sigma - \dim \tau}(1 - t)^{\codim \sigma}
 \prod_{\rho_{i} \subseteq \sigma \backslash \tau} (1 - t)/(1 - s^{\lambda(b_{i})}t).
\end{equation}
The weighted $h$-vector is a rational function in $\mathbb{Q}(s^{1/N'},t)$, for some positive integer $N'$.
Note that if $\lambda \equiv 0$ or if we set $s = 1$, then $h_{\tau}^{\lambda}(s,t)$ is equal to the usual $h$-vector $h_{\tau}(t)$ of $\Sigma_{\tau}$.
By expanding and collecting terms, we have the following equality, 
\begin{equation}\label{snitch}
t^{\codim \tau} h_{\tau}^{\lambda}(s^{-1},t^{-1}) =
\sum_{\tau \subseteq \sigma} (t - 1)^{\codim \sigma} \prod_{\rho_{i} \subseteq \sigma \backslash \tau} (t - 1)/(s^{\lambda(b_{i})}t - 1).
\end{equation}

\begin{lemma}\label{LemmaA}
Suppose that $\tau$ is a cone not contained in the boundary of $\Sigma$. Then
\begin{equation*}
h_{\tau}^{\lambda}(s,t) = t^{\codim \tau} h_{\tau}^{\lambda}(s^{-1},t^{-1}).
\end{equation*}
\end{lemma}
\begin{proof}
This is an application of M\"obius inversion (see, for example, \cite{StaEnumerative}).
Let $\mathcal{P}$ be the poset consisting of the cones in $\Sigma_{\tau}$ and of a maximal element $\{ \Sigma_{\tau} \}$.
By Lemma \ref{Lefty}, $h_{\tau}(t) = t^{\codim \tau}h_{\tau}(t^{-1})$. Substituting $t  = 0$ gives
$\sum_{\tau \subseteq \sigma} (- 1)^{\codim \sigma} = 1$. It follows that we can compute the M\"obius function of $\mathcal{P}$ inductively.
M\"obius inversion says that if $f: \mathcal{P} \rightarrow A$ is a function, for some abelian group $A$, then
\begin{equation}\label{eqnM}
f(\{ \Sigma_{\tau} \}) = g(\{\Sigma_{\tau}\} ) + \sum_{\tau \subseteq \sigma} (-1)^{\codim \sigma + 1}g(\sigma),
\end{equation}
where, for each $p$ in $\mathcal{P}$, $g(p) = \sum_{q \leq p} f(q)$. 

We set $f(\{ \Sigma_{\tau}\}) = 0$ and
\begin{equation*}
f(\sigma) = (t- 1)^{\codim \sigma} \prod_{\rho_{i} \subseteq \sigma \backslash \tau} (t - 1)/(s^{\lambda(b_{i})}t - 1).
\end{equation*}
Then   $g(\{\Sigma_{\tau}\} ) = t^{\codim \tau} h_{\tau}^{\lambda}(s^{-1},t^{-1})$ by (\ref{snitch}) and we calculate
\begin{align*}
g(\sigma) &= \sum_{\tau \subseteq \sigma' \subseteq \sigma} f(\sigma') \\
&= (t - 1)^{\codim \tau}\sum_{\tau \subseteq \sigma' \subseteq \sigma}  \prod_{\rho_{i} \subseteq \sigma' \backslash \tau} 1/(s^{\lambda(b_{i})}t - 1)  \\
&= (t - 1)^{\codim \tau} \prod_{\rho_{i} \subseteq \sigma \backslash \tau} (1 + 1/(s^{\lambda(b_{i})}t - 1)  ) \\
&= s^{\sum_{\rho_{i} \in \sigma \backslash \tau}\lambda(b_{i})} t^{\dim \sigma - \dim \tau}f(\sigma).
\end{align*}
By (\ref{pretty}) and  (\ref{eqnM}),
\[
h_{\tau}^{\lambda}(s,t)  = \sum_{\tau \subseteq \sigma} (-1)^{\codim \sigma} g(\sigma)
= g(\{\Sigma_{\tau}\} ) = t^{\codim \tau} h_{\tau}^{\lambda}(s^{-1},t^{-1}).
\]
\end{proof}

As in \cite{BCSOrbifold}, for each non-zero cone $\tau$ of $\Sigma$, set
\begin{equation}\label{mcbox}
\BOX(\btau) = \{ v \in N \mid  v = \sum_{\rho_{i} \subseteq \tau} q_{i}b_{i} \textrm{  for some  }
0 < q_{i} < 1 \}.
\end{equation}
We set $\BOX(\mbox{\boldmath$\{ 0 \}$}) = \{0\}$ and $\BOX(\bSigma) = \cup_{\tau \in \Sigma} \BOX(\btau)$. Given any $v$ in $|\Sigma| \cap N$, let $\sigma(v)$ be the cone in
$\Sigma$ containing $v$ in its relative interior. As in \cite[p.7]{MPEhrhart}, $v$ has a unique decomposition
\begin{equation}\label{eqnD}
v = \{ v \} + v' + \sum_{\rho_{i} \subseteq \sigma(v) \backslash \tau} b_{i},
\end{equation}
where $\{v \}$ lies in $\BOX(\btau)$ for some $\tau \subseteq \sigma(v)$ and $v'$ is a linear combination of the $\{ b_{i} \, | \, \rho_{i} \subseteq \sigma(v) \}$
with  non-negative integer coefficients. We think of $\{ v \}$ as the `fractional part' of $v$.
Using this decomposition, we compute a local formula for the weighted $\delta$-vector of $Q$. The method of proof is the same as that of Theorem 1.3 in \cite{MPEhrhart} and Theorem 1.2 in
\cite{PayEhrhart}. In fact, Proposition \ref{PropA} can be deduced from Theorem 1.2 in \cite{PayEhrhart}.

\begin{proposition}\label{PropA}
The weighted $\delta$-vector of $Q$ is a rational function in $\mathbb{Q}(s^{1/N},t)$, for some positive integer $N$, and
has the form
\begin{equation*}
\delta^{\lambda}(s,t)  = \sum_{\tau \in \Sigma} h_{\tau}^{\lambda}(s,t) \sum_{ v \in \BOX(\btau)} s^{w_{\lambda}(v)}
t^{\lceil \psi(v) \rceil} \prod_{\rho_{i} \subseteq \tau} (t -1)/(s^{\lambda(b_{i})}t -1) .
\end{equation*}
\end{proposition}
\begin{proof}
\begin{align*}
\delta^{\lambda}(s,t)  
&= (1 -t)^{d + 1} \sum_{m \geq 0} \sum_{ v \in mQ \cap N }s^{w_{\lambda}(v)}t^{m} \\
&= (1 -t)^{d + 1}  \sum_{ v \in |\Sigma| \cap N   } \sum_{m \geq \psi(v)}s^{w_{\lambda}(v)}t^{m} \\
&= (1 - t)^{d} \sum_{ v \in |\Sigma| \cap N}  s^{w_{\lambda}(v)}t^{\lceil \psi(v) \rceil}.
\end{align*}
Now consider the decomposition (\ref{eqnD}). We have 
\begin{equation*}
w_{\lambda}(v) = w_{\lambda}(\{ v \} ) + w_{\lambda}(v') + \sum_{\rho_{i} \subseteq \sigma(v) \backslash \tau} \lambda(b_{i}).
\end{equation*}
Also
\begin{equation*}
 \lceil \psi(v) \rceil =  \lceil \psi(\{ v \}) \rceil + \psi(v') + \dim \sigma(v) - \dim \tau.
\end{equation*}
We obtain the following expression for $\delta^{\lambda}(t)$
\begin{equation*}
(1 - t)^{d} \sum_{ \substack{ \tau \in \Sigma \\ v \in \BOX(\btau)} } s^{w_{\lambda}(v)}t^{\lceil \psi(v) \rceil} \sum_{\tau \subseteq \sigma}
s^{\sum_{\rho_{i} \subseteq \sigma \backslash \tau} \lambda(b_{i})} t^{\dim \sigma - \dim \tau} \prod_{ \rho_{i} \subseteq \sigma} 1/( 1 - s^{\lambda(b_{i})}t).
\end{equation*}
Substituting in (\ref{pretty}) and rearranging gives the result.
\end{proof}


\begin{exa}\label{exo}
If $\lambda \equiv 0$, then 
\begin{equation*}
\delta^{0}(s,t)  = \sum_{\substack{\tau \in \Sigma \\ v \in \BOX(\btau)}} s^{\psi(v) - \lceil \psi(v) \rceil}t^{\lceil \psi(v) \rceil}
h_{\tau}(t),
\end{equation*}
where $h_{\tau}(t)$ is the usual $h$-vector of $\Sigma_{\tau}$. In this case, $\delta^{0}(s,t) \in \mathbb{Z}[s^{1/N},s^{-1/N},t]$, for some positive integer $N$.
By Lemma \ref{Lefty}, the coefficients of
$h_{\tau}(t)$ are non-negative integers. Hence the coefficients of $\delta^{0}(s,t)$ are non-negative integers.
Recall that, for any rational
number $-1 < k \leq 0$, we can recover the polynomial $\delta^{0}_{k}(t)$ as the coefficient of $s^{k}$ in $\delta^{0}(s,t)$.
Note that the degree of $h_{\tau}(t)$ is at most $\codim \tau$ and $\lceil \psi(v) \rceil \leq \dim \tau$ for any $v$ in $\BOX(\btau)$.
It follows that $\delta^{0}_{k}(t)$ is a polynomial in $\mathbb{Z}[t]$ of degree less than or equal to $d$,
with non-negative coefficients.
\end{exa}

\begin{exa}\label{exaA}
Recall that we can recover the Ehrhart $\delta$-vector of $Q$ as $\delta_{Q}(t) = \delta^{\lambda}(1,t)$.
Then Proposition \ref{PropA} gives the local formula of Betke and McMullen for $\delta_{Q}(t)$ 
(Theorem 1 in \cite{BMLattice}),           
\begin{equation*}
\delta_{Q}(t)  = \sum_{ \substack{ \tau \in \Sigma \\ v \in \BOX(\btau)}} t^{\lceil \psi(v) \rceil}
h_{\tau}(t).
\end{equation*}
Note that the origin in $N$ corresponds to  a contribution of $h_{0}(t)$ in the above sum.
It follows that $\delta_{Q}(t)$ is a polynomial of degree less than or equal to $d$ with non-negative coefficients and constant term $1$ \cite{BMLattice}. If $\Sigma$ is complete,
we conclude from Lemma \ref{Lefty}
that $\delta_{Q}(t)$ is a polynomial of degree $d$ with positive integer coefficients \cite{BMLattice}.
\end{exa}


\begin{cor}\label{Coral}
Suppose that $\lambda$ satisfies the additional condition
\begin{equation}\label{eqnn}
\lambda(v) \geq 0,
\end{equation}
for every $v$ in $|\Sigma|$.
Then for every rational number $k$, $\delta^{\lambda}_{k}(t)$ is a polynomial in $\mathbb{Z}[t]$.
\end{cor}
\begin{proof}
We know that $\delta^{\lambda}_{k}(t)$ is the coefficient of $s^{k}$ in $\delta^{\lambda}(s,t)$ and is a power series in $\mathbb{Z}[[t]]$.
We will show that it has bounded degree. By Proposition \ref{PropA},
\begin{equation*}
\delta^{\lambda}(s,t)  = \sum_{\tau \in \Sigma} h_{\tau}^{\lambda}(s,t)
\sum_{ v \in \BOX(\btau)} s^{w_{\lambda}(v)}t^{\lceil \psi(v) \rceil} \prod_{\rho_{i} \subseteq \tau} (t -1)/(s^{\lambda(b_{i})}t -1) .
\end{equation*}
Expanding the right hand side gives
\begin{equation*}
(1 - t)^{d}\sum_{ \substack{ \tau \in \Sigma \\ v \in \BOX(\btau)} } s^{w_{\lambda}(v)}t^{\lceil \psi(v) \rceil} \sum_{\tau \subseteq \sigma}
s^{\sum_{\rho_{i} \subseteq \sigma \backslash \tau} \lambda(b_{i})} t^{\dim \sigma - \dim \tau} \prod_{ \rho_{i} \subseteq \sigma} 1/( 1 - s^{\lambda(b_{i})}t).
\end{equation*}
If a monomial $t^{l}s^{k}$ appears in the expansion of this expression, then $k$ must have the form
\begin{equation}\label{cafe}
k = w_{\lambda}(v) + \sum_{\rho_{i} \subseteq \sigma \backslash \tau} \lambda(b_{i}) + \sum_{i = 1}^{r} \alpha_{i} \lambda(b_{i}),
\end{equation}
for some $v$ in $\BOX(\bSigma)$ and $\alpha_{i}$ non-negative integers that are equal to zero if $\lambda(b_{i}) = 0$, and such that
\begin{equation*}
l \leq \lceil \psi(v) \rceil + 2d + \sum_{i = 1}^{r}\alpha_{i}.
\end{equation*}
It follows from Condition (\ref{eqnn}) that for a fixed $k$, there are only finitely many possibilities for $\alpha_{i}$ such that (\ref{cafe}) holds. Therefore
$l$ is bounded.
\end{proof}

By a standard argument, this is equivalent
to the following corollary.
\begin{cor}\label{coast}
Suppose that $\lambda$ satisfies the additional condition
\begin{equation*}
\lambda(v) \geq 0,
\end{equation*}
for every $v$ in $|\Sigma|$.
For every rational number $k$ and for every $m$ sufficiently large (depending on $k$),
$f^{\lambda}_{k}(m)$ is a polynomial function in $m$, of degree less than or equal to $d$. 
\end{cor}
\begin{proof}
Fix a rational number $k$. By Corollary \ref{Coral}, we can write
\begin{equation*}
F^{\lambda}_{k}(t) = P_{k}(t)/(1 - t)^{d + 1} + Q_{k}(t),
\end{equation*}
where $P_{k}(t) = p_{0,k} + p_{1,k}t + \cdots + p_{d,k}t^{d}$ and $Q_{k}(t)$ are polynomials.
Expanding the right hand side gives
\begin{equation*}
F^{\lambda}_{k}(t) = \sum_{j = 0}^{d} p_{j,k} \sum_{m \geq 0} \binom{m + d}{d} t^{m + j} + Q_{k}(t).
\end{equation*}
Hence, for $m > \deg Q_{k}(t)$,
\begin{equation*}
f^{\lambda}_{k}(m) = \sum_{j = 0}^{d} p_{j,k} \binom{m + d - j}{d}.
\end{equation*}
\end{proof}

\begin{exa}\label{extra}
As in Example \ref{exo}, suppose that $\lambda \equiv 0$. We have seen that for every rational number $k$,
$\delta^{0}_{k}(t)$ is a polynomial in $\mathbb{Z}[t]$ with non-negative coefficients, of degree less than or equal to $d$. We claim that
$f^{0}_{k}(m)$ is either identically zero or a polynomial of degree $d$ in $\mathbb{Q}[t]$ with positive leading coefficient. This follows from the above proof, which
shows that $f^{0}_{k}(m)$ is a polynomial of degree less than or equal to $d$, and that the coefficient of $t^{d}$ is $\sum_{j =0}^{d}p_{j,k}/d!$,
where the $p_{j,k}$ are the non-negative coefficients of
$\delta^{0}_{k}(t)$.
\end{exa}

Suppose that $\lambda$ satisfies the additional condition
\begin{equation}\label{leche}
\lambda(b_{i}) > -1 \textrm{ for } i = 1, \ldots, r.
\end{equation}
In this case we define
\begin{equation*}
h_{\tau}^{\lambda}(t) := h_{\tau}^{\lambda}(t,t),
\end{equation*}
\begin{equation*}
\delta^{\lambda}(t) := \delta^{\lambda}(t,t).
\end{equation*}
By Proposition \ref{PropA}, we have the following expression for $\delta^{\lambda}(t)$,
\begin{equation}\label{oak}
\delta^{\lambda}(t)  = \sum_{\tau \in \Sigma} h_{\tau}^{\lambda}(t) \sum_{ v \in \BOX(\btau)} t^{\psi(v) + \lambda(v)}
\prod_{\rho_{i} \subseteq \tau} (t -1)/(t^{\lambda(b_{i}) + 1} -1) .
\end{equation}
It follows from (\ref{leche}) and (\ref{pretty}) that there is a positive integer $N'$ such that $h_{\tau}^{\lambda}(t)$ lies in $\mathbb{Z}[[t^{1/N'}]]$ for all $\tau$ in $\Sigma$.
By  (\ref{leche}), $w_{\lambda}(v) + \lceil \psi(v) \rceil =  \psi(v) + \lambda(v)$ is non-negative for all $v$ in $|\Sigma| \cap N$.
Hence  
(\ref{oak}) implies that  $\delta^{\lambda}(t) \in \mathbb{Z}[[t^{1/N}]]$, for some positive integer $N$. In this case, we will often abuse notation and call $\delta^{\lambda}(t)$ the
\emph{weighted $\delta$-vector} associated to $\bSigma$ and $\lambda$.

We now describe a well-known involution $\iota$ on $|\Sigma| \cap N$. Consider a cone $\tau$ in $\Sigma$ and $v$ in $\BOX(\btau)$. Then $v$ can be uniquely written in the form
$v = \sum_{\rho_{i} \subseteq \tau} q_{i}b_{i}$, for some $0< q_{i} <1$. We define
\begin{equation}\label{invo}
\iota = \iota_{\bSigma} : \BOX(\btau) \rightarrow \BOX(\btau)
\end{equation}
\begin{equation*}
\iota(v) = \sum_{\rho_{i} \subseteq \tau} (1 - q_{i})b_{i}.
\end{equation*}
As in (\ref{eqnD}), every $v$ in $|\Sigma| \cap N$ can be uniquely written in the form
$v = \{ v \} + \tilde{v}$, where $\{ v \}$ is in $\BOX(\btau)$ for some $\tau \subseteq \sigma(v)$ and $\tilde{v}$ is in $N_{\sigma(v)}$. Here 
$\sigma(v)$ is the cone of $\Sigma$ containing $v$ in
its relative interior.
Then $\iota$ extends to an involution on $|\Sigma| \cap N$
\begin{equation*}
\iota = \iota_{\bSigma} : |\Sigma| \cap N \rightarrow |\Sigma| \cap N
\end{equation*}
\begin{equation*}
\iota(v) = \iota(\{ v \}) + \tilde{v}.
\end{equation*}
Using (\ref{oak}) and after rearranging and collecting terms, we can write $t^{d}\delta^{\lambda}(t^{-1})$ as
\begin{equation*}  
\sum_{  \substack{\tau \in \Sigma \\ v \in \BOX(\btau)} } t^{\codim \tau}h_{\tau}^{\lambda}(t^{-1})
t^{\sum_{\rho_{i} \subseteq \tau}\lambda(b_{i}) + \dim \tau - \psi(v) -   \lambda(v)}
\prod_{\rho_{i} \subseteq \tau} (t -1)/(t^{\lambda(b_{i}) + 1} -1). 
\end{equation*}  
Note that for any $v$ in $\BOX(\btau)$,
\[
\sum_{\rho_{i} \subseteq \tau}\lambda(b_{i}) + \dim \tau - \psi(v) -   \lambda(v) = \psi(\iota(v)) + \lambda(\iota(v)),
\]
where $\iota$ is the involution (\ref{invo}). Hence 
\begin{equation}\label{googly}
t^{d}\delta^{\lambda}(t^{-1}) = \sum_{\tau \in \Sigma} t^{\codim \tau}h_{\tau}^{\lambda}(t^{-1}) \sum_{ v \in \BOX(\btau)}
t^{\psi(v) + \lambda(v)}
\prod_{\rho_{i} \subseteq \tau} (t -1)/(t^{\lambda(b_{i}) + 1} -1).
\end{equation}

\begin{cor}\label{beautiful}
Assume that $\Sigma$ is a complete fan and $\lambda(b_{i}) > -1$ for $i = 1, \ldots, r$.
Then
\begin{equation*}
\delta^{\lambda}(t) = t^{d} \delta^{\lambda}(t^{-1}).
\end{equation*}
\end{cor}
\begin{proof}
By Lemma \ref{LemmaA}, for any cone $\tau$ in $\Sigma$, $h^{\lambda}_{\tau}(t) = t^{\codim \tau}h^{\lambda}_{\tau}(t^{-1})$.
The result follows by comparing expressions (\ref{oak}) and (\ref{googly}).
\end{proof}

We have the following change of variables formula for weighted $\delta$-vectors.
A geometric proof involving motivic integration is given in  \cite{YoWeightII}.   
We write $\psi = \psi_{\Sigma}$ and $\delta^{\lambda}(t) = \delta^{\lambda}_{\Sigma}(t)$.


\begin{proposition}\label{change}
Let $N$ be a lattice of rank $d$. Let $\bSigma = (N, \Sigma, \{ b_{i} \})$ and $\btriangle = (N , \triangle, \{ b'_{j} \})$ be stacky fans such that $|\Sigma| = |\triangle|$. Let $\lambda$ be
a piecewise $\mathbb{Q}$-linear function with respect to $\Sigma$ satisfying $\lambda(b_{i}) > -1$ for every $b_{i}$, and set
\[
\lambda' = \lambda +  \psi_{\Sigma} - \psi_{\triangle} .
\]
If $\lambda'$ is piecewise $\mathbb{Q}$-linear with respect to $\triangle$ and satisfies $\lambda'(b'_{j}) > -1$ for every $b'_{j}$,
then
$\delta^{\lambda}_{\Sigma}(t) = \delta^{\lambda'}_{\triangle}(t)$.
\end{proposition}
\begin{proof}
By letting $s = t$ in the first calculation in the proof of Proposition \ref{PropA}, we see that
\[
\delta^{\lambda}_{\Sigma}(t)  = (1 - t)^{d} \sum_{ v \in |\Sigma| \cap N}  t^{\psi_{\Sigma}(v) + \lambda(v)},
\]
which gives,
\[
\delta^{\lambda'}_{\triangle}(t)  = (1 - t)^{d} \sum_{ v \in |\triangle| \cap N}  t^{\psi_{\triangle}(v) + \lambda'(v)} = \delta^{\lambda}_{\Sigma}(t).
\]
\end{proof}

\section{Weighted Ehrhart Reciprocity}\label{Ehr}

The goal of this section is to investigate the case when $\lambda \equiv 0$. In this case,
$\delta^{0}(t)$ is a polynomial of degree at $d$ with rational powers and non-negative integer coefficients,
and we can recover the Ehrhart $\delta$-vector $\delta_{Q}(t)$ from $\delta^{0}(t)$.
While $\delta_{Q}(t)$  is not symmetric in its coefficients, Corollary \ref{beautiful} implies if $\Sigma$ is complete,
then $\delta^{0}(t) = t^{d}\delta^{0}(t^{-1})$. The main point is that we can exploit this symmetry  to deduce facts about the Ehrhart $\delta$-vector.



Note that the weight function $w_{0}(v) = \psi(v) - \lceil \psi(v) \rceil$ takes values between $-1$ and $0$.
By Example \ref{exo}, for each rational number $-1< k \leq 0$, we can write
\begin{equation*}
\delta^{0}_{k}(t) = \delta_{d,k}t^{d} + \delta_{d - 1,k}t^{d - 1} + \cdots + \delta_{0,k},
\end{equation*}
for some non-negative integers $\delta_{i,k}$.
Since the Ehrhart $\delta$-vector decomposes as $\delta_{Q}(t) = \sum_{k \in (-1,0]} \delta^{0}_{k}(t)$, we have, with the notation of  (\ref{otation}),
\begin{equation}\label{bling}
\delta_{i} = \sum_{k \in (-1,0]} \delta_{i,k},
\end{equation}

Throughout this section we will set $s = t$, so that the weighted $\delta$-vector is given by
\begin{equation*}
\delta^{0}(t) = \sum_{k \in (-1,0]} \delta^{0}_{k}(t) t^{k} = \sum_{k \in (-1,0]} \sum_{i = 0}^{d} \delta_{i,k} t^{i + k}.
\end{equation*}
By (\ref{bling}), $\delta_{i}$ is the sum of the coefficients of
$t^{j}$ in $\delta^{0}(t)$ for $i - 1 < j \leq i$.
By Example \ref{exo},
\begin{equation}\label{equat}
\delta^{0}(t)  =  \sum_{\substack{\tau \in \Sigma \\ v \in \BOX(\btau)}}  t^{\psi(v)}
h_{\tau}(t).
\end{equation}
Later we will see that the coefficients of $\delta^{0}(t)$
are dimensions of  orbifold Chow groups of a toric stack (Theorem \ref{boo}).

\begin{rem}
The weight function $w_{0}$ and hence the weighted $\delta$-vector $\delta^{0}(t)$ are determined by the underlying space of the lattice complex
$Q$.
Recall, from Remark \ref{impo}, that any lattice polytope $P$, after translation, has the form $P= Q_{\bSigma}$, for some stacky fan $\bSigma$.
In this case, $\delta^{0}(t)$ is determined by $P$ and the choice of a lattice point $\alpha$ in $P$.
\end{rem}

\begin{rem}
Note that the non-zero lattice points of weight $0$ are those that lie in a facet of $\partial (mQ)$ not containing $0$, for some positive integer $m$.
Let $\partial Q_{0}$ denote the union of the facets in $\partial Q$ not containing the origin.
Consider the lattice $N \times \mathbb{Z}$ and the lattice complex
\[
K_{0} = \{ (v, \mu) \in (N \times \mathbb{Z})_{\mathbb{R}} \mid 0 < \mu \leq 1, \, v \in \partial(\mu Q)_{0} \} \cup \{ 0 \}.
\]
We can interpret $f^{0}_{0}(m)$ as the number of lattice points in $mK_{0}$ and $\delta^{0}_{0}(t)$ as the $\delta$-vector associated to $K_{0}$  (c.f. \cite{YoInequalities}).
\end{rem}

\begin{exa}
Let $P$ be a lattice $d$-simplex containing the origin in its interior. Let $\Sigma$ be the fan over the faces of the boundary of $P$, with the appropriate choice of $\{b_{i}\}$.
One can verify that $h_{\tau}(t) = 1 + t + \cdots + t^{\codim \tau}$ for any cone $\tau$ in $\Sigma$. Hence
\[
\delta^{0}(t)  = \sum_{\substack{\tau \in \Sigma \\ v \in \BOX(\btau)}} t^{\psi(v)} ( 1 + t + \cdots + t^{\codim \tau} ) .
\]
\end{exa}



\begin{cor}\label{cozzie}
Suppose $\Sigma$ is a complete fan.
Then $\delta^{0}_{0}(t)$ has degree $d$ and positive integer coefficients. 
For $-1 < k < 0$, we can write
$\delta^{0}_{k}(t) = t\tilde{\delta}_{k}(t)$.
Then
\begin{equation*}
\delta^{0}(t) = t^{d}\delta^{0}(t^{-1})
\end{equation*}
\begin{equation*}
\delta^{0}_{0}(t) = t^{d}\delta^{0}_{0}(t^{-1})
\end{equation*}
\begin{equation*}
\tilde{\delta}_{k}(t) = t^{d - 1}\tilde{\delta}_{-1 - k}(t^{-1}),
\end{equation*}
for $-1 < k < 0$.
\end{cor}
\begin{proof}
By considering the contribution of $0 \in \BOX \{ 0 \}$ in (\ref{equat}) and using Lemma \ref{Lefty},
we see that $\delta^{0}_{0}(t)$ has degree $d$ and positive integer coefficients.
By (\ref{equat}), $\delta^{0}_{k}(t)$ has no constant term for $-1 < k < 0$.
By Corollary \ref{beautiful}, we have $\delta^{0}(t) = t^{d} \delta^{0}(t^{-1})$.
We can write
\begin{equation*}
\delta^{0}(t) = \sum_{k \in (-1,0]} \delta^{0}_{k}(t) t^{k}.
\end{equation*}
\begin{align*}
t^{d}\delta^{0}(t^{-1}) &= \sum_{k \in (-1,0]} t^{d}\delta^{0}_{k}(t^{-1}) t^{-k} \\
&= t^{d}\delta^{0}_{0}(t^{-1}) + \sum_{k \in (-1,0)} t^{d + 1}\delta^{0}_{k}(t^{-1}) t^{-1-k} \\
&= t^{d}\delta^{0}_{0}(t^{-1}) + \sum_{k \in (-1,0)} t^{d + 1}\delta^{0}_{-1 - k}(t^{-1}) t^{k}. \\
\end{align*}
Comparing the expressions above yields
\begin{equation*}
\delta^{0}_{0}(t) = t^{d}\delta^{0}_{0}(t^{-1})
\end{equation*}
\begin{equation*}
\delta^{0}_{k}(t) = t^{d + 1}\delta^{0}_{-1-k}(t^{-1})
\end{equation*}
for $-1 < k < 0$. Since each $\delta^{0}_{k}(t)$ is a polynomial of degree less than
or equal to $d$, the corollary follows.
\end{proof}

\begin{rem}
Suppose $\Sigma$ is complete and set
\begin{equation*}
a(t) = \delta^{0}_{0}(t), \; b(t) = \sum_{k \in (-1,0)} \tilde{\delta}_{k}(t).
\end{equation*}
Since $\delta_{Q}(t) = \sum_{-1 < k \leq 0} \delta^{0}_{k}(t)$, we have a decomposition
\begin{equation*}
\delta_{Q}(t) = a(t) + tb(t).
\end{equation*}
This decomposition is due to  Betke and McMullen  (Theorem 5 \cite{BMLattice}).
We give an analogous result in the case when $\Sigma$ is not necessarily complete in \cite{YoInequalities}.
\end{rem}

\begin{rem}
We can exploit the symmetry properties of the above decomposition to translate inequalities
between the coefficients of $\delta_{Q}(t)$ into inequalities between the coefficients of $\delta^{0}_{0}(t)$. This is explained in \cite{YoInequalities} in a more general setting.
\end{rem}

We give a reformulation of Corollary \ref{cozzie}. We note that the $k = 0$ case can be deduced from Ehrhart Reciprocity (see Remark \ref{Geelong}).

\begin{theorem}[Weighted Ehrhart Reciprocity]\label{propB}
Suppose $\Sigma$ is a complete fan.
For every rational number $- 1 < k \leq 0$, $f^{0}_{k}(m)$ is either identically zero or
a polynomial of degree $d$ in
$\mathbb{Q}[t]$ with positive leading coefficient.
For any positive integer $m$,
\begin{displaymath}
f^{0}_{k}(-m) = \left\{ \begin{array}{ll}
(-1)^{d}f^{0}_{k}(m - 1) & \textrm{ if } k = 0 \\
(-1)^{d}f^{0}_{-1-k}(m) & \textrm{ if } -1 < k < 0.
\end{array} \right.
\end{displaymath}
\end{theorem}
\begin{proof}
Fix a rational number $k$. The first statement follows from Example \ref{extra}.  We see from the proof of
Corollary \ref{coast} that
\begin{equation*}
f^{0}_{k}(m) = \sum_{j = 0}^{d} \delta_{j,k} \binom{m + d - j}{d}.
\end{equation*}
Hence
\begin{equation*}
f^{0}_{k}(-m) = (-1)^{d}\sum_{j = 0}^{d} \delta_{j,k} \binom{m + j - 1}{d} \textrm{ for } m \geq 1.
\end{equation*}
We define
\begin{equation*}
F_{k}(t) = \sum_{m \geq 0} f^{0}_{k}(m)t^{m}, \tilde{F}_{k}(t) = \sum_{m \geq 1} f^{0}_{k}(-m)t^{m}.
\end{equation*}
We compute
\begin{align*}
-F_{k}(t^{-1}) &= -\delta^{0}_{k}(t^{-1})/(1 - t^{-1})^{d + 1} \\
&= (-1)^{d}t^{d + 1}\delta^{0}_{k}(t^{-1})/(1-t)^{d + 1} \\
&= (-1)^{d} \sum_{j = 0}^{d} \sum_{m \geq 0} \delta_{j,k} \binom{m + d}{d} t^{m + d + 1 - j} \\
&= \sum_{j = 0}^{d} \sum_{m' \geq d +1 - j}  (-1)^{d}\delta_{j,k} \binom{m' + j -1}{d} t^{m'}\\
&= \sum_{m' \geq 1} \sum_{j = 0}^{d}   (-1)^{d}\delta_{j,k} \binom{m' + j -1}{d} t^{m'}\\
&= \tilde{F}_{k}(t).
\end{align*}
It follows from Corollary \ref{cozzie} (after substituting in definitions) that
\begin{align*}
tF_{0}(t) &= (-1)^{d + 1}  F_{0}(t^{-1}) \\
F_{k}(t) &= (-1)^{d + 1} F_{-1 -k}(t^{-1}) \textrm{ for } -1 < k < 0.
\end{align*}
The result follows by comparing coefficients of the expressions
\begin{align*}
\tilde{F}_{0}(t) &= (-1)^{d} t F_{0}(t) \\
\tilde{F}_{k}(t) &= (-1)^{d} F_{-1-k}(t) \textrm{ for } -1 < k < 0.
\end{align*}
\end{proof}

The above result should be viewed as the weighted version of Ehrhart Reciprocity.
In particular, we verify below that Ehrhart Reciprocity is a consequence.
We emphasise that this proof 
can be interpreted as a reformulation of previous proofs of Ehrhart's theorem.


\begin{cor}[Ehrhart Reciprocity \cite{EhrDemonstration}]\label{Ehrt}

For every positive integer $m$, \\
$(-1)^{d}f_{Q}(-m)$ is equal to the number of lattice points in the interior of $mQ$.
\end{cor}
\begin{proof}
It follows from the definitions that
\begin{equation}\label{rubbish}
(-1)^{d}f_{Q}(-m) = \sum_{k \in (-1,0] \cap \, \mathbb{Q}} (-1)^{d}f^{0}_{k}(-m)
\end{equation}
\begin{equation}\label{rubbish2}
|\Int(mQ) \cap N| = f^{0}_{0}(m - 1) + \sum_{k \in (-1,0) \cap \, \mathbb{Q}} f^{0}_{k}(m).
\end{equation}

The fan $\Sigma$ induces a lattice triangulation $\mathcal{T}$ of $\partial Q$.
There exists a positive integer $n$ and a lattice point $\alpha$ in the interior of $nQ$ such that, after translating $\alpha$ to the origin,
the collection of cones over the faces of $\mathcal{T}$ form a simplicial fan $\Sigma'$. It follows from (\ref{rubbish}), (\ref{rubbish2}) and Example \ref{exo} that
the functions $|\Int(mQ) \cap N|$ and $(-1)^{d}f_{Q}(-m)$ are polynomials in $m$.
Hence, after replacing $Q$ by a multiple and replacing $\Sigma$ by
$\Sigma'$, we may assume that $\Sigma$ is complete. 
The result now follows by applying Theorem \ref{propB} to (\ref{rubbish}) in order to obtain (\ref{rubbish2}).

\end{proof}

\begin{rem}\label{Geelong}
This result was conjectured by Ehrhart in 1959 and proved by him in \cite{EhrDemonstration}. Another proof was given by MacDonald in \cite{MacPolynomials}.
If $\Sigma$ is complete, then
\begin{equation}\label{Jolly}
f_{Q}(m) - (-1)^{d}f_{Q}(-m) = f^{0}_{0}(m) - f^{0}_{0}(m -1),
\end{equation}
since both sides are equal to the number of lattice points in $\partial(mQ)$. With the notation of (\ref{poly}), it follows that $c_{d - 1}$ is half the surface area of
$Q$, normalised with respect to the sublattice on each facet of $Q$. Note that $c_{d}$ is the normalised volume of $Q$ in $N$. These facts were established in \cite{EhrLinearII} and \cite{MacVolume}.
Using (\ref{Jolly}) and applying a short induction,  Ehrhart Reciprocity implies that $f^{0}_{0}(-m) =
(-1)^{d}f^{0}_{0}(m - 1)$ for any positive integer $m$, which was proved in Theorem \ref{propB} (c.f. \cite{YoInequalities}).
\end{rem}

\begin{rem}\label{TMac}
In a similar way, one can show that Weighted Ehrhart Reciprocity implies Ehrhart Reciprocity for rational polytopes (see, for example,
\cite{StaEnumerative}). More specifically, fix a positive integer $r$ and set $Q' = (1/r)Q$. For $0 \leq l < r$, one verifies from the definitions that
\[
f_{Q'}(l + mr) = \sum_{l/r -1 < k \leq 0} f_{k}(m) + \sum_{-1 < k \leq l/r - 1} f_{k}(m + 1).
\]
It follows that $f_{Q'}(m)$ is a quasipolynomial with period dividing $r$ \cite{StaEnumerative}. As in the proof of Corollary
\ref{Ehrt}, after replacing $Q$ by $sQ$ for some positive integer $s$ coprime to $r$, we can apply Theorem \ref{propB} to deduce that, for any positive integer $m$,
$(-1)^{d}f_{Q'}(-m)$
is equal to the number of lattice points in the interior of $mQ'$.

Let $P$ be a rational polytope  and fix a positive integer $r$ such that $rP$ is a lattice polytope. The function $f_{P}(m)$ is a quasipolynomial with period dividing $r$ \cite{StaEnumerative}. After replacing $P$ by $sP$ for some positive integer $s$ coprime to $r$ and translating by a lattice point,
we may assume that $P$ contains the origin. Setting $Q' = P$ and $Q = rP$, we deduce that
$(-1)^{d}f_{P}(-m)$
is equal to the number of lattice points in the interior of $mP$.


\end{rem}

The following result of Hibi was originally shown to be a consequence of Ehrhart Reciprocity \cite{HibDual}. In a similar way, it follows from Theorem \ref{propB}.
Recall that $\psi$ is a piecewise $\mathbb{Q}$-linear function.

\begin{cor}[\cite{HibDual}]\label{gogo}
If $\Sigma$ is complete,
then $\delta_{Q}(t) = t^{d}\delta_{Q}(t^{-1})$ if and only if
$\psi$ is a piecewise linear function.
\end{cor}
\begin{proof}
Observe that $\psi$ is piecewise linear if and only if
$w_{0} \equiv 0$.
If $w_{0} \equiv 0$, then by Corollary \ref{cozzie}, $\delta_{Q}(t) = \delta^{0}(t) = t^{d}\delta^{0}(t^{-1}) =   t^{d}\delta_{Q}(t^{-1})$.
Conversely,
suppose that $\delta_{Q}(t) = t^{d}\delta_{Q}(t^{-1})$.  Assume $w_{0}$ is not identically $0$. Choose $j + 1$ minimal such that $\delta_{j + 1,k} > 0$ for some $-1 < k < 0$.
By Corollary \ref{cozzie} and (\ref{bling}),
\begin{equation*}
\delta_{j} = \delta_{j,0} + \sum_{k \in (-1,0)} \delta_{j,k} =  \delta_{j,0}.
\end{equation*}
\begin{equation*}
\delta_{d - j} = \delta_{d - j,0} + \sum_{k \in (-1,0)} \delta_{d- j,k} =  \delta_{j,0} + \sum_{k \in (-1,0)} \delta_{j + 1,k} > \delta_{j}.
\end{equation*}
This is a contradiction.
\end{proof}

\begin{rem}
If $P$ is a lattice polytope then $\delta_{0}$ is $1$ and $\delta_{d}$ is the number of interior lattice points of $P$.
A lattice polytope $P$ is \emph{reflexive} if it contains the origin in its interior and $\psi$ is piecewise linear.
Then Corollary \ref{gogo} says that $\delta_{P}(t) = t^{d}\delta_{P}(t^{-1})$ if and only if $P$ is the translate of a reflexive polytope \cite{HibDual}.
\end{rem}

We conclude this section  by proving a lemma which allows us to compute examples when $d = 2$ and $\Sigma$ is complete, and then presenting a corresponding example.
Recall that $\partial Q_{0}$ denotes the union of the facets of $\partial Q$ not containing the origin.

\begin{lemma}\label{Kansas}
The weighted $\delta$-vector $\delta^{0}(t)$ is a polynomial of degree less than or equal to $d$ with rational powers
and non-negative integer coefficients. If $\Sigma$ is complete, then  $\delta^{0}(t) = t^{d} \delta^{0}(t^{-1})$.
\begin{enumerate}
\item $\delta^{0}(1) = d! \vol_{d}(Q)$.
\item The constant coefficient in $\delta^{0}(t)$ is 1.
\item  For $0 < l < 1$, the coefficient of $t^{l}$ in $\delta^{0}(t)$ is $| \{ v \in Q \cap N \mid \psi(v) = l\}|$.
\item The coefficient of $t$ in $\delta^{0}(t)$ is $| \partial Q_{0} \cap N| - d$.
\end{enumerate}
\end{lemma}
\begin{proof}
We established the initial claims in Corollary \ref{coast} and Corollary \ref{cozzie}.
We showed in the proof of Corollary \ref{coast} that $f^{0}_{k}(m)$ is a polynomial of degree $d$ with leading coefficient
$\sum_{j = 0}^{d} \delta_{j,k}/d!$. Hence
\begin{equation*}
d!\vol_{d}(Q) = \sum_{k \in (-1,0]} \sum_{j = 0}^{d} \delta_{j,k} = \delta^{0}(1).
\end{equation*}

For the other claims we compare both sides of the expression
\begin{equation*}
\delta^{0}(t) = (1 - t)^{d + 1} \sum_{m \geq 0} \sum_{v \in mQ \cap N} t^{w_{0}(v) + m}.
\end{equation*}
The constant coefficient on the right hand side is $1$. For $0 < l < 1$, the coefficient on the right hand side is $| \{ v \in Q \cap N \mid w_{0}(v) + 1 = l \}|$.
Note that if $v$ lies in $Q \cap N$ and
$w_{0}(v) \neq 0$ then $\lceil \psi(v) \rceil = 1$. Finally, the coefficient of $t$ on the right hand side is $|\{ v \in Q \cap N \mid w_{0}(v) = 0 \}| - (d+ 1)$.
The only elements of $Q \cap N$ of weight zero
are the origin and the elements of $\partial Q_{0} \cap N$.
\end{proof}

\begin{rem}\label{cuando}
When $d = 2$ and $\Sigma$ is complete,
one can show that
\[ f^{0}_{0}(m) = |\partial Q \cap N|m(m + 1)/2 + 1 \] \[ f^{0}_{k}(m) = (f^{0}_{k}(1) +  f^{0}_{-1 - k}(1)) m^{2}/2 + ( f^{0}_{k}(1) -  f^{0}_{-1 - k}(1)) m/2 \textrm{ for } k \neq 0. \]
\end{rem}

\begin{exa}

Let $N = \mathbb{Z}^{2}$ and let $\Sigma$ be the complete fan with primitive integer vectors
$(1,0)$,$(1,3)$,
$(0,1)$, $(-2,3)$,$(-2,1)$,$(-1,0)$ and $(0,-1)$, and set $a_{i}$ to be $1,1,2,1,1,2$ and $1$ respectively.
Since $\Sigma$ is a complete fan, weighted Ehrhart Reciprocity holds and the weighted $\delta$-vector $\delta^{0}(t)$ is symmetric.
The example is illustrated in the diagram below, in which we have marked the
lattice points of non-zero weight in $2Q$. The computations were made using
Lemma \ref{Kansas} and Remark \ref{cuando}. Observe that removing the ray through $(-2,1)$ does not affect the results below.
\[ \delta^{0}(t) = t^{2} + 3t^{3/2} + t^{5/4} + 8t + t^{3/4} + 3t^{1/2} + 1 \]
\[ \delta_{Q}(t) = 5t^{2} + 12t + 1 \]

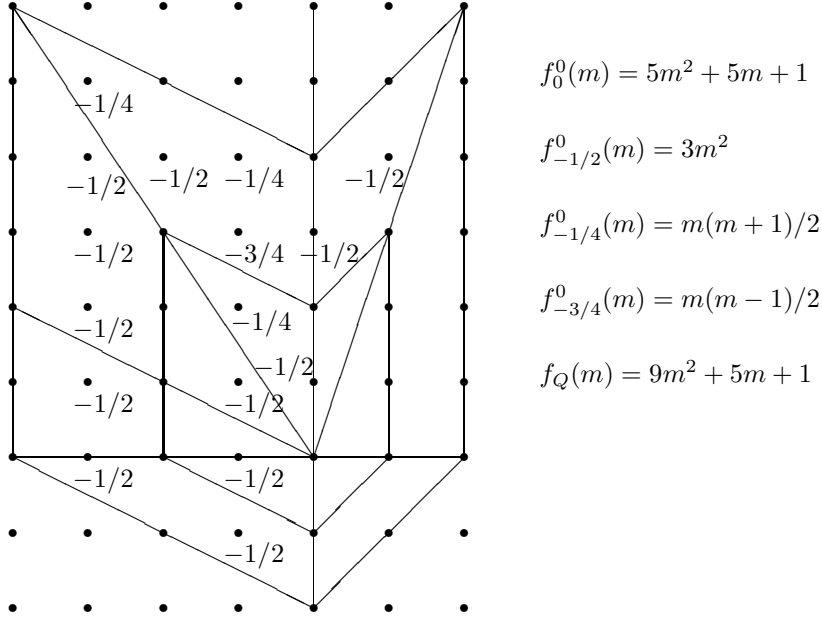
\setlength{\unitlength}{1cm}
\begin{picture}(12,8.5)

\put(7,7){$f^{0}_{0}(m) = 5m^{2} + 5m + 1$}
\put(7,6){$f^{0}_{-1/2}(m) = 3m^{2}$}
\put(7,5){$f^{0}_{-1/4}(m) = m(m + 1)/2$}
\put(7,4){$f^{0}_{-3/4}(m) = m(m - 1)/2$}
\put(7,3){$f_{Q}(m) = 9m^{2} + 5m + 1$}

\put(0,0){\circle*{0.1}}
\put(0,1){\circle*{0.1}}
\put(0,2){\circle*{0.1}}
\put(0,3){\circle*{0.1}}
\put(0,4){\circle*{0.1}}
\put(0,5){\circle*{0.1}}
\put(0,6){\circle*{0.1}}
\put(0,7){\circle*{0.1}}
\put(0,8){\circle*{0.1}}
\put(1,0){\circle*{0.1}}
\put(1,1){\circle*{0.1}}
\put(1,2){\circle*{0.1}}
\put(1,3){\circle*{0.1}}
\put(1,4){\circle*{0.1}}
\put(1,5){\circle*{0.1}}
\put(1,6){\circle*{0.1}}
\put(1,7){\circle*{0.1}}
\put(1,8){\circle*{0.1}}
\put(2,0){\circle*{0.1}}
\put(2,1){\circle*{0.1}}
\put(2,2){\circle*{0.1}}
\put(2,3){\circle*{0.1}}
\put(2,4){\circle*{0.1}}
\put(2,5){\circle*{0.1}}
\put(2,6){\circle*{0.1}}
\put(2,7){\circle*{0.1}}
\put(2,8){\circle*{0.1}}
\put(3,0){\circle*{0.1}}
\put(3,1){\circle*{0.1}}
\put(3,2){\circle*{0.1}}
\put(3,3){\circle*{0.1}}
\put(3,4){\circle*{0.1}}
\put(3,5){\circle*{0.1}}
\put(3,6){\circle*{0.1}}
\put(3,7){\circle*{0.1}}
\put(3,8){\circle*{0.1}}
\put(4,0){\circle*{0.1}}
\put(4,1){\circle*{0.1}}
\put(4,2){\circle*{0.1}}
\put(4,3){\circle*{0.1}}
\put(4,4){\circle*{0.1}}
\put(4,5){\circle*{0.1}}
\put(4,6){\circle*{0.1}}
\put(4,7){\circle*{0.1}}
\put(4,8){\circle*{0.1}}
\put(5,0){\circle*{0.1}}
\put(5,1){\circle*{0.1}}
\put(5,2){\circle*{0.1}}
\put(5,3){\circle*{0.1}}
\put(5,4){\circle*{0.1}}
\put(5,5){\circle*{0.1}}
\put(5,6){\circle*{0.1}}
\put(5,7){\circle*{0.1}}
\put(5,8){\circle*{0.1}}
\put(6,0){\circle*{0.1}}
\put(6,1){\circle*{0.1}}
\put(6,2){\circle*{0.1}}
\put(6,3){\circle*{0.1}}
\put(6,4){\circle*{0.1}}
\put(6,5){\circle*{0.1}}
\put(6,6){\circle*{0.1}}
\put(6,7){\circle*{0.1}}
\put(6,8){\circle*{0.1}}

\put(3.2,3.1){$-1/2$}
\put(3.8,4.6){$-1/2$}
\put(4.4,5.6){$-1/2$}
\put(2.8,1.6){$-1/2$}
\put(0.8,1.6){$-1/2$}
\put(0.8,2.6){$-1/2$}
\put(0.8,3.6){$-1/2$}
\put(0.8,4.6){$-1/2$}
\put(0.7,5.5){$-1/2$}
\put(2.8,2.6){$-1/2$}
\put(2.8,4.6){$-3/4$}
\put(2.8,5.6){$-1/4$}
\put(1.8,5.6){$-1/2$}
\put(0.8,6.6){$-1/4$}
\put(2.9,3.7){$-1/4$}
\put(2.8,0.6){$-1/2$}


\linethickness{0.075mm}
\put(4,0){\line(0,1){8}}
\put(0,2){\line(1,0){6}}
\put(4,2){\line(-2,3){4}}
\put(4,2){\line(-2,1){4}}
\put(4,2){\line(1,3){2}}

\linethickness{0.2mm}
\put(2,5){\line(2,-1){2}}
\put(2,5){\line(0,-1){3}}
\put(2,2){\line(0,1){1}}
\put(4,4){\line(1,1){1}}
\put(5,2){\line(-1,-1){1}}
\put(2,2){\line(2,-1){2}}
\put(6,2){\line(-1,-1){2}}
\put(0,2){\line(2,-1){4}}
\put(0,8){\line(0,-1){6}}
\put(5,5){\line(0,-1){3}}
\put(6,2){\line(0,1){6}}
\put(0,8){\line(2,-1){4}}
\put(4,6){\line(1,1){2}}

\end{picture}

\end{exa}

\section{Orbifold Cohomology}\label{orbifold}

We will now prove our geometric interpretation of the coefficients of the Ehrhart $\delta$-vector.
Recall that $\Sigma$ is a simplicial, $d$-dimensional fan with convex support. In this case, we have a combinatorial description of the Betti numbers of the corresponding
toric variety $X = X(\Sigma)$.

\begin{lemma}\label{AFL}
Let $\Sigma$ is a simplicial, $d$-dimensional fan with convex support.
Let $X = X(\Sigma)$ be the toric variety associated to $\Sigma$. Then $X$ has no odd cohomology over $\mathbb{Q}$ and
$\dim  H^{2i}( X, \mathbb{Q})$ is equal to the coefficient of $t^{i}$ in the $h$-vector $h_{\Sigma}(t)$ of $\Sigma$.
If $\Sigma$ is complete, then
$\dim  H^{2i}( X, \mathbb{Q}) > 0$ for  $i = 0, \ldots, d$.
\end{lemma}
\begin{proof}
The case when $\Sigma$ is complete is due to Danilov (Theorem 10.8 \cite{DanGeometry}). 
Suppose $\Sigma$ is not complete. Let $\rho$ be a ray in the interior of $-|\Sigma|$. Let
$\triangle$ be the fan with cones given by the cones of $\Sigma$ as well as the cones generated by $\rho$ and a face of $\Sigma$ contained in $\partial |\Sigma|$.
Consider the complete, simplicial toric variety $Y = Y(\triangle)$.
Let $D = D(\triangle_{\rho})$ be the $\mathbb{Q}$-Cartier torus-invariant divisor corresponding to $\rho$. Then $D$ is simplicial and complete and $Y \smallsetminus D = X$.
By considering the long exact sequence of cohomology with compact supports, 
we have a diagram,
\[
\xymatrix{   &  & A^{i}(Y, \mathbb{Q}) \ar[r]^{\alpha} \ar[d] & A^{i}(D, \mathbb{Q}) \ar[r] \ar[d] & 0 & \\
0 \ar[r] & H_{\com}^{2i}( X, \mathbb{Q})  \ar[r] &  H_{\com}^{2i}( Y, \mathbb{Q})  \ar[r] &  H_{\com}^{2i}( D, \mathbb{Q}) \ar[r] & H_{\com}^{2i + 1}( X, \mathbb{Q}) \ar[r] & 0.}
\]
Since a complete toric variety has no odd cohomology, the bottom row is exact.
The vertical maps take a cycle to its corresponding cohomology class and both maps are isomorphisms
by Theorem 10.8 of \cite{DanGeometry}.
The map $\alpha$ `restricts' cycles of $Y$ to cycles on the $\mathbb{Q}$-Cartier divisor $D$ (p33 \cite{FulIntersection}).
Let $\gamma$ be a cone in $\partial |\Sigma|$ corresponding to a $T$-invariant subvariety $V(\gamma)$ of $Y$.  If we set $\sigma = \gamma + \rho$, then $\sigma$ corresponds to a $T$-invariant subvariety $V_{\rho}(\sigma)$ of $D$ and the set-theoretic intersection of $V(\gamma)$ and $D$ is $V_{\rho}(\sigma)$. It follows that  $\alpha([V(\gamma)])$ is
a positive multiple of $[V_{\rho}(\sigma)]$ \cite{FulIntroduction}.
By Proposition 10.3 of \cite{DanGeometry}, $A^{*}(D, \mathbb{Q})$ is generated by classes of the form $[V_{\rho}(\sigma)]$ and hence $\alpha$ is surjective and the diagram above has exact rows.
We conclude that
$H_{\com}^{2i + 1}( X, \mathbb{Q}) = 0$ and $\dim H_{\com}^{2i}( X, \mathbb{Q})$ is equal to the $i^{\textrm{th}}$ coefficient of
\[
t^{d}h_{\triangle}(t^{-1}) - t^{d}h_{\triangle_{\rho}}(t^{-1}) = \sum_{\tau \in \Sigma} (t - 1)^{\codim \tau} = t^{d}h_{\Sigma}(t^{-1}).
\]
By Poincar\'e duality, $\dim H^{2i}( X, \mathbb{Q})$  is equal to the coefficient of $t^{i}$ in $h_{\Sigma}(t)$.
\end{proof}

We can associate to the stacky fan $\bSigma = (N, \Sigma, \{ b_{i} \} )$ a Deligne-Mumford toric stack $\mathcal{X} = \mathcal{X}(\bSigma)$ over $\mathbb{C}$ with
coarse moduli space $X = X(\Sigma)$  \cite{BCSOrbifold}. Let $\mathcal{Y}_{1}, \ldots, \mathcal{Y}_{t}$ denote the connected components of the corresponding inertia stack
$\mathcal{I} \mathcal{X}$. There is a degree shifting function
\[
\iota_{\mathcal{X}}: \mathcal{I} \mathcal{X} \rightarrow \mathbb{Q},
\]
which is constant on connected components. Let $\overline{\mathcal{Y}}_{i}$ denote the coarse moduli space of $\mathcal{Y}_{i}$.
For $i \in \mathbb{Q}$, Chen and Ruan \cite{CRNew} defined the $i^{\textrm{th}}$ orbifold cohomology group of $\mathcal{X}$ by
\begin{equation*}
H_{\orb}^{i}(\mathcal{X}, \mathbb{Q}) = \bigoplus_{j = 1}^{t} H^{i - 2\iota_{\mathcal{X}}(\mathcal{Y}_{j}   )}(\overline{\mathcal{Y}}_{j}, \mathbb{Q}).
\end{equation*}
Similarly, we can consider the orbifold cohomology $H_{\orb, \com}^{*}(\mathcal{X}, \mathbb{Q})$ of $\mathcal{X}$ with compact support.
Chen and Ruan established Poincar\'e duality between $H_{\orb}^{*}(\mathcal{X}, \mathbb{Q})$ and $H_{\orb, \com}^{*}(\mathcal{X}, \mathbb{Q})$ (Proposition 3.3.1 \cite{CRNew}).

Borisov, Chen and Smith (Proposition 4.7 \cite{BCSOrbifold}) show that the connected components of $\mathcal{I} \mathcal{X}$ are indexed by the elements of  $\BOX(\bSigma)$.
Moreover, if $v$ in $\BOX(\btau)$ corresponds to the connected component $\mathcal{Y}_{v}$, then $\iota_{\mathcal{X}}( \mathcal{Y}_{v} ) = \psi(v)$ and
$\overline{\mathcal{Y}}_{v} = X(\Sigma_{\tau})$. Hence
\begin{equation}\label{nearly}
H_{\orb}^{2i}(\mathcal{X}, \mathbb{Q}) = \bigoplus_{\tau \in \Sigma} \bigoplus_{v \in \BOX(\btau)} H^{2(i - \psi(v))}( X(\Sigma_{\tau}), \mathbb{Q}).
\end{equation}

\begin{rem}
Similarly, Abramovich, Graber and Vistoli \cite{AGVAlgebraic} defined the  \emph{orbifold Chow ring} $A_{\orb}^{*}(\mathcal{X}, \mathbb{Q})$ of $\mathcal{X}$. For $i \in \mathbb{Q}$,
\[
A_{\orb}^{i}(\mathcal{X}, \mathbb{Q}) =  \bigoplus_{j = 1}^{t} H^{i - \iota_{\mathcal{X}}(\mathcal{Y}_{j}   )}(\overline{\mathcal{Y}}_{j}, \mathbb{Q}).
\]
The cohomology ring with rational coefficients and Chow ring with rational coefficients of a simplicial, complete toric variety are isomorphic
(Theorem 10.8 \cite{DanGeometry}). Hence, if $\Sigma$ is complete, then $A_{\orb}^{i}(\mathcal{X}, \mathbb{Q}) \cong H_{\orb, \com}^{2i}(\mathcal{X}, \mathbb{Q})$.
\end{rem}

We finally arrive at the main result of this section. The case when $w_{0} \equiv 0$ is proved in \cite{MPEhrhart}.

\begin{theorem}\label{boo}
The coefficient of $t^{j}$ in $\delta^{0}(t)$ is equal to $\dim_{\mathbb{Q}} H_{\orb}^{2j}(\mathcal{X}(\bSigma), \mathbb{Q})$. Moreover, the coefficient $\delta_{i}$ of $t^{i}$ in the $\delta$-vector
$\delta_{Q}(t)$ is a sum of dimensions of orbifold cohomology groups,
\begin{equation*}
\delta_{i}  = \sum_{2i - 2< j \leq 2i} \dim_{\mathbb{Q}} H_{\orb}^{j}(\mathcal{X}(\bSigma), \mathbb{Q}).
\end{equation*}
\end{theorem}
\begin{proof}
By (\ref{equat}),
\[
\delta^{0}(t) = \sum_{\tau \in \Sigma} \sum_{ v \in \BOX(\btau)} h_{\tau}(t) t^{\psi(v)}.
\]
The first assertion now follows from (\ref{nearly}) and Lemma \ref{AFL}.
The second statement then follows from
(\ref{bling}).
\end{proof}

\begin{rem}
By Poincar\'e duality, the coefficient of $t^{j}$ in $t^{d}\delta^{0}(t^{-1})$
 is equal to the dimension of the $(2j)^{\textrm{th}}$ orbifold cohomology group of $\mathcal{X}$ with compact support.
\end{rem}

\begin{rem}\label{bias}
When $\Sigma$ is complete, we showed in  Corollary \ref{cozzie} that 
\begin{equation*}
\delta^{0}(t) = t^{d}\delta^{0}(t^{-1}).
\end{equation*}
By Theorem \ref{boo}, we can interpret this symmetry as a consequence of Poincar\'e duality for orbifold cohomology. 
In particular, since Weighted Ehrhart Reciprocity (Theorem \ref{propB}) is equivalent to Corollary \ref{cozzie}, this provides a geometric proof of Weighted Ehrhart Reciprocity.
\end{rem}

A corollary of Theorem \ref{boo} is the following result which interprets the coefficients of the Ehrhart $\delta$-vector of a lattice polytope as
dimensions of orbifold cohomology groups of a $(d + 1)$-dimensional orbifold. More specifically,
let $P$ be a $d$-dimensional lattice polytope in $N$ and
fix a lattice triangulation $\mathcal{T}$ of $P$.  If
$\sigma$ denotes the cone over $P \times \{ 1 \}$ in $(N \times \mathbb{Z})_{\mathbb{R}}$, then $\mathcal{T}$ determines a simplicial fan refinement $\triangle$ of $\sigma$. The corresponding toric variety $Y = Y(\triangle)$ has a canonical stack structure: if
$w_{1}, \ldots, w_{s}$ are the primitive integer vectors of the rays of $\triangle$, then the corresponding stacky fan is
$(N \times \mathbb{Z}, \triangle, \{w_{i}\})$. We will write  $H_{\orb}^{2i}(Y, \mathbb{Q})$ for the $2i^{\textrm{th}}$ orbifold cohomology group of the canonical stack associate to $Y$.

\begin{theorem}\label{chance}
Let $P$ be a $d$-dimensional lattice polytope and let $\mathcal{T}$ be a lattice triangulation of $P$ corresponding to a $(d + 1)$-dimensional toric variety $Y$ as above. The Ehrhart $\delta$-vector of $P$ has the form
\[
\delta_{P}(t) = \sum_{i = 0}^{d} \dim_{\mathbb{Q}} H_{\orb}^{2i}(Y, \mathbb{Q}) t^{i}.
\]
\end{theorem}
\begin{proof}
With the notation of the previous discussion, let $\mathcal{Y}$ denote the toric stack associated to the stacky fan
$(N \times \mathbb{Z}, \triangle, \{w_{i}\})$. In this case, $\psi: |\triangle| \rightarrow \mathbb{R}$ is the restriction of the projection $N_{\mathbb{R}} \times \mathbb{R} \rightarrow \mathbb{R}$ to $|\triangle|$, and hence $Q = \{ v \in |\triangle| \mid \psi(v) \leq 1 \}$ is the convex hull of $P \times \{1\}$ and the origin, called the \emph{pyramid} over $P$.
Since the weight function $w_{0}(v) = \psi(v) - \lceil \psi(v) \rceil$ is identically zero on $|\triangle| \cap (N \times \mathbb{Z})$,
the weighted $\delta$-vector $\delta^{0}(t)$ is just the usual $\delta$-vector $\delta_{Q}(t)$. It is a standard fact that $\delta_{P}(t) = \delta_{Q}(t)$ (see, for example, Remark 2.6 \cite{BatLattice}). On the other hand, Theorem \ref{boo} implies that
$\delta_{P}(t) = \delta_{Q}(t)  = \delta^{0}(t) = \sum_{i = 0}^{d} \dim_{\mathbb{Q}} H_{\orb}^{2i}(\mathcal{Y}, \mathbb{Q}) t^{i}$. 
\end{proof}

\begin{rem}
If $\mathcal{T}$ is a unimodular triangulation of $P$, then $Y$ is smooth and $H_{\orb}^{2i}(Y, \mathbb{Q}) = H^{2i}(Y, \mathbb{Q})$. In this case, the above theorem and Lemma \ref{AFL} imply the well-known fact that $\delta_{P}(t)$ is equal to the $h$-vector of $\mathcal{T}$ \cite{HibEhrhart}.
\end{rem}

\section{A Toric Proof of Weighted Ehrhart Reciprocity}\label{toric}

We have provided a combinatorial proof of Weighted Ehrhart Reciprocity (Theorem \ref{propB}) as well as a geometric proof via orbifold cohomology (Remark \ref{bias}). In this section, we give a third proof in the case when $Q$ is a lattice polytope. More specifically, we show that Weighted Ehrhart Reciprocity can be deduced from Serre Duality as well as some vanishing theorems for ample divisors on toric varieties due to Musta\c t\v a \cite{MusVanishing}. This proof generalises the toric proof of Ehrhart Reciprocity given in Section 4.4 of \cite{FulIntroduction}.

Throughout this section we will assume that $\Sigma$ is a complete fan and
\begin{equation*}
P := Q = \{ v \in N_{\mathbb{R}} \mid \psi(v) \leq 1 \}
\end{equation*}
is a lattice polytope (containing the origin in its interior). By definition, 
for every rational number $- 1 < k \leq 0$
and for any positive integer $m$,
\begin{equation}\label{savings}
f^{0}_{k}(m) - f^{0}_{k}(m - 1) = | \partial (m + k)P \cap N |.
\end{equation}
Also, $f^{0}_{0}(0) = 1$ and $f^{0}_{k}(0) = 0$ for $-1 < k < 0$.
Our goal is to provide a toric proof of the following version of Weighted Ehrhart Reciprocity.
\begin{theorem}\label{India}
Let $P$ be a $d$-dimensional lattice polytope containing the origin in its interior.
For every rational number $- 1 < k \leq 0$,
$f^{0}_{k}(m)$ is 
a polynomial in $m$ of degree at most $d$
and for any positive integer $m$,
\begin{displaymath}
f^{0}_{k}(-m) = \left\{ \begin{array}{ll}
(-1)^{d}f^{0}_{-1-k}(m) & \textrm{ if } -1 < k < 0 \\
(-1)^{d}f^{0}_{k}(m - 1) & \textrm{ if } k = 0.
\end{array} \right.
\end{displaymath}
\end{theorem}

We first recall some facts about toric varieties and refer the reader to \cite{FulIntroduction} for the relevant details.
A $d$-dimensional lattice polytope $P$ in $N$, containing the origin in its interior, determines a $d$-dimensional, projective
toric variety $Y$ over $\mathbb{C}$ and an effective ample
torus-invariant divisor $D$ on $Y$. If $M = \Hom (N, \mathbb{Z})$ denotes the dual lattice to $N$, then
the normal fan to $P$ in $M_{\mathbb{R}}$ determines the toric variety $Y$.
Let $u_{1}, \ldots, u_{s}$ denote the primitive integer vectors along the rays of the normal fan, corresponding to the
torus-invariant prime divisors $D_{1}, \ldots, D_{s}$ of $Y$. If we write $D = \sum_{i =1 }^{s} a_{i}D_{i}$, then
$a_{i} = -\min_{v \in P \cap N} \langle u_{i}, v \rangle \in \mathbb{Z}_{> 0}$.
Given a torus-invariant $\mathbb{Q}$-divisor $E = \sum_{i = 1}^{s} b_{i} D_{i}$, consider the (possibly empty) polytope
\[
P_{E} := \{ v \in N_{\mathbb{R}} \mid \langle u_{i} , v \rangle + b_{i} \geq 0 \textrm{ for } i = 1, \ldots, s \}.
\]
It can be verified that $\lambda P = P_{ \lambda D }$ for any rational number $\lambda > 0$.
Every lattice point $v$ in $N$ corresponds to a character $\chi^{v}$ on the torus contained in $Y$. In particular, we may view
$\chi^{v}$ as a rational function on $Y$. If we let $\lfloor E \rfloor = \sum_{i = 1}^{s} \lfloor b_{i} \rfloor D_{i}$ denote the round down of $E$,
then the global sections of $E$ are given by
\begin{equation}\label{coffee}
H^{0}(Y, \mathcal{O}(\lfloor E \rfloor))  = \bigoplus_{v \in P_{E} \cap N} \mathbb{C} \chi^{v}.
\end{equation}
In particular, $\dim H^{0}(Y, \mathcal{O}(\lfloor E \rfloor)) = | P_{E} \cap N|$.
We identify the canonical divisor of $Y$ with
$K_{Y} = - \sum_{i = 1}^{s} D_{i}$ and write $\lceil E \rceil = \sum_{i = 1}^{s} \lceil b_{i} \rceil D_{i}$ for the round up of $E$.

\begin{lemma}\label{letter}
If $E = \sum_{i = 1}^{s} b_{i} D_{i}$ is a torus-invariant $\mathbb{Q}$-divisor on $Y$ then
\[
H^{0}(Y, \mathcal{O}(K_{Y} + \lceil E \rceil))  = \bigoplus_{v \in \Int P_{E} \cap N} \mathbb{C} \chi^{v}.
\]
In particular, for each $-1 < k \leq 0$ and for every positive integer $m$,
\begin{equation*}
f^{0}_{k}(m) - f^{0}_{k}(m - 1) = \dim H^{0}(Y, \mathcal{O}(\lfloor (m + k)D \rfloor)) -
\dim H^{0}(Y, \mathcal{O}(K_{Y} + \lceil (m + k)D \rceil)).
\end{equation*}
\end{lemma}
\begin{proof}
Observe that $v$ in $N$ lies in the interior of $P_{E}$ if and only if
\[
\langle u_{i} , v \rangle + b_{i} > 0 \textrm{ for } i = 1, \ldots, s.
\]
This condition holds if and only if
\[
\langle u_{i} , v \rangle + \lceil b_{i} \rceil - 1 \geq 0 \textrm{ for } i = 1, \ldots, s.
\]
We conclude that $\Int P_{E} \cap N = P_{K_{Y} + \lceil E \rceil} \cap N$ and the first statement follows from (\ref{coffee}).
By (\ref{savings}), for each $-1 < k \leq 0$ and for every positive integer $m$,
\begin{align*}
f^{0}_{k}&(m) - f^{0}_{k}(m - 1) = | (m + k)P \cap N | - |\Int(m + k)P \cap N | \\
&= \dim H^{0}(Y, \mathcal{O}(\lfloor (m + k)D \rfloor)) -
\dim H^{0}(Y, \mathcal{O}(K_{Y} + \lceil (m + k)D \rceil)).
\end{align*}
\end{proof}

We recall the following vanishing theorem due to Musta\c t\v a. A $\mathbb{Q}$-divisor $E$ on $Y$ is ample if $mE$ is an ample divisor for some positive integer $m$.

\begin{theorem}[Corollary 2.5 \cite{MusVanishing}]\label{lion}
Let $Y$ be a projective toric variety and let $E$ be an ample torus-invariant $\mathbb{Q}$-divisor on $Y$. For $i > 0$,
\begin{enumerate}
\item
$H^{i}(Y, \mathcal{O}(K_{Y} + \lceil E \rceil)) = 0$ (Kawamata-Viehweg vanishing)
\item
$H^{i}(Y, \mathcal{O}(\lfloor E \rfloor)) = 0$.
\end{enumerate}
\end{theorem}

For any divisor $D'$ on $Y$, let $\chi(Y, D') = \sum_{i \geq 0} (-1)^{i} \dim H^{i}(Y, \mathcal{O}(D'))$ denote the Euler characteristic
of $D'$.
By the above vanishing theorem and Lemma \ref{letter}, for any positive integer $m$,
\begin{equation}\label{queue}
f^{0}_{k}(m) - f^{0}_{k}(m - 1) = \chi(Y, \mathcal{O}(\lfloor (m + k)D \rfloor)) -
\chi(Y, \mathcal{O}(K_{Y} + \lceil (m + k)D \rceil)).
\end{equation}

The following fact is due to Snapper.

\begin{theorem}[\cite{SnaMultiples}]
Let $X$ be a complete variety of dimension $d$ over an algebraically closed field $k$. If $\mathcal{F}$ is a coherent sheaf on
$X$ and $\mathcal{L}$ is a line bundle on $X$,  then there is a
polynomial $Q(t)$ of degree
at most $d$ such that $Q(m) = \chi(X, \mathcal{F} \otimes \mathcal{L}^{m})$ for every integer $m$. 
\end{theorem}

If we apply the above result to the coherent sheaf $\mathcal{O}(\lfloor kD \rfloor)$ and the line bundle $\mathcal{O}(D)$ on $Y$, we deduce that there is a polynomial $Q_{1}(t)$ of degree at most $d$ such
that $Q_{1}(m) = \chi(Y, \mathcal{O}(\lfloor (m + k)D \rfloor))$ for each integer $m$. Moreover,
the coefficient of $t^{d}$ in
$Q_{1}(t)$ is the intersection number $\rk (\mathcal{O}(\lfloor kD \rfloor)) \cdot D^{d}/d! = D^{d}/d! > 0$ \cite{KleToward}.
Similarly, there is a polynomial $Q_{2}(t)$ of degree $d$ with leading term $D^{d}/d!$ such that
$Q_{2}(m) = \chi(Y, \mathcal{O}(K_{Y} + \lceil (m + k)D \rceil))$ for each integer $m$.
By (\ref{queue}), for any positive integer $m$, $f^{0}_{k}(m) - f^{0}_{k}(m - 1) = Q_{1}(m) - Q_{2}(m)$. It now follows from standard
arguments (see, for example, p. 49 \cite{HarAlgebraic}) that $f^{0}_{k}(m)$ is a polynomial in $m$ of degree at most $d$.

We have the following application of Serre Duality. 

\begin{lemma}\label{tarheels}
If $S_{k}(m) = f^{0}_{k}(m) - f^{0}_{k}(m - 1)$ then 
\begin{displaymath}
(-1)^{d + 1}S_{k}(-m) = \left\{ \begin{array}{ll}
S_{-1-k}(m + 1) & \textrm{ if } -1 < k < 0 \\
S_{k}(m) & \textrm{ if } k = 0.
\end{array} \right.
\end{displaymath}
\end{lemma}
\begin{proof}
Since both sides of (\ref{queue}) are polynomials in $m$, $(-1)^{d + 1}S_{k}(-m)$ is equal to
\[
-(-1)^{d}\chi(Y, \mathcal{O}(\lfloor (-m + k)D \rfloor)) +
(-1)^{d}\chi(Y, \mathcal{O}(K_{Y} + \lceil (-m + k)D \rceil)).
\]
By Serre Duality (see, for example, Corollary 3.7.7 \cite{HarAlgebraic}), 
this is equal to
\[
-\chi(Y, \mathcal{O}(K_{Y} - \lfloor (-m + k)D \rfloor)) +
\chi(Y, \mathcal{O}(-\lceil (-m + k)D \rceil)).
\]
Since for any real number $a$, $-\lfloor -a \rfloor = \lceil a \rceil$,
\[
(-1)^{d + 1}S_{k}(-m) = - \chi(Y, \mathcal{O}(K_{Y} + \lceil (m - k)D \rceil)) +
\chi(Y, \mathcal{O}(\lfloor (m - k)D \rfloor)).
\]
When $k = 0$, we get $(-1)^{d + 1}S_{0}(-m) = S_{0}(m)$. For $-1 < k \leq 0$, after writing $m - k = m + 1 + (-1 - k)$, we see
that $(-1)^{d + 1}S_{k}(-m) = S_{-1-k}(m + 1)$.
\end{proof}

We will now complete our proof of Theorem \ref{India} by induction on $m$.

\begin{proof}
We have seen that for each $-1 <  k \leq 0$, $f^{0}_{k}(m)$ is a polynomial in $m$ of degree at most $d$.  We will first prove Theorem  \ref{India} in the case when $m = 1$.
Recall that $f^{0}_{0}(0) = 1$ and $f^{0}_{k}(0) = 0$ for $-1 < k < 0$.
By (\ref{queue}) and Serre Duality,
\[
f^{0}_{0}(0) - f^{0}_{0}(-1) = \chi(Y, \mathcal{O}_{Y}) -
\chi(Y, \mathcal{O}(K_{Y})) = (1 - (-1)^{d})\chi(Y, \mathcal{O}_{Y}).
\]
It follows from Theorem \ref{lion} that $\chi(Y, \mathcal{O}_{Y}) = 1$ and we conclude that
$f^{0}_{0}(-1) = (-1)^{d} = (-1)^{d}f^{0}_{0}(0)$ as desired. When $-1 < k < 0$, Lemma \ref{tarheels} implies
that $(-1)^{d + 1}S_{k}(0) = S_{-1-k}(1)$. That is,
\[
(-1)^{d + 1}(f^{0}_{k}(0) - f^{0}_{k}(-1)) = f^{0}_{-1-k}(1) - f^{0}_{-1-k}(0),
\]
and hence $(-1)^{d}f^{0}_{k}(-1) = f^{0}_{-1-k}(1)$. This completes the proof when $m = 1$.

Now consider the case when $m > 1$. By Lemma \ref{tarheels} and induction on $m$,
\begin{align*}
(-1)^{d}f^{0}_{0}(-m) &= (-1)^{d + 1}S_{0}(-m + 1) + (-1)^{d}f^{0}_{0}(-(m - 1)) \\
&= S_{0}(m - 1) + f^{0}_{0}(m - 2) \\
&= f^{0}_{0}(m - 1).
\end{align*}
Similarly, when $-1< k < 0$,
\begin{align*}
(-1)^{d}f^{0}_{k}(-m) &= (-1)^{d + 1}S_{k}(-m + 1) + (-1)^{d}f^{0}_{k}(-(m - 1)) \\
&= S_{-1- k}(m) + f^{0}_{-1-k}(m - 1) \\
&= f^{0}_{-1-k}(m).
\end{align*}

\end{proof}

\bibliographystyle{amsplain}
\bibliography{alan}

\end{document}